\documentclass[11pt]{amsart}
\UseRawInputEncoding
\usepackage{amsmath,amssymb,amsfonts,url,mathptmx}
\numberwithin{equation}{section}

\newtheorem{thm}{Theorem}[section]
\newtheorem{lem}{Lemma}[section]
\newtheorem{rem}{Remark}[section]

\begin{document}
\title[Singular Liouville Equation]{Estimates for Liouville equation with quantized singularities}

\author{Juncheng Wei}
\address{Department of Mathematics \\ University of British Columbia\\ Vancouver, BC V6T1Z2, Canada} \email{jcwei@math.ubc.ca }

\author{Lei Zhang} \footnote{The research of J. Wei is partially supported by NSERC of Canada. Lei Zhang is partially supported by a Simons Foundation Collaboration Grant}
\address{Department of Mathematics\\
        University of Florida\\
        1400 Stadium Rd\\
        Gainesville FL 32611}
\email{leizhang@ufl.edu}

\date{\today}

%%%%%%%%%%%%%%%%%%%%%%%%%%%%%%%%%%%%%%%%%%%%%
\begin{abstract} For Liouville equations with singular sources, the interpretation of the equation and its impact are most significant if the singular sources are quantized: the strength of each Dirac mass is a mutliple of $4\pi$. However the study of bubbling solutions around a quantized singular source is particularly challenging: near the singular source, the spherical Harnack inequality may not hold and there are multiple local maximums of bubbling solutions all swarming to the singular source. In this article we seek to provide a complete understanding of the blowup picture in this core difficulty and we establish two major types of results: first we prove that not only the first derivatives of coefficient functions tend to zero  at the singular source, the second derivatives also have a vanishing estimate. Second we derive pointwise estimates for bubbling solutions to be approximated by global solutions, which have crucial applications in a number of important projects. Since bubbling solutions near a singular source are very commonly observed in geometry and physics, it seems that the estimates in this article can also be applied to many related equations and systems with various backgrounds.
\end{abstract}

%%%%%%%%%%%%%%%%%%%%%%%%%%%%%%%%%%%%%%%%%%%%%
\maketitle

\section{Introduction}
The main purpose of this article is to study bubbling solutions of the following singular Liouville equation
\begin{equation}\label{main-eq}
\Delta u+ h(x) e^{u}=4\pi \alpha \delta_0\quad \mbox{ in }\quad \Omega\subset \mathbb R^2
\end{equation}
where $\Omega$ is an open, bounded subset of $\mathbb R^2$ that contains the origin, $\alpha>-1$ is a constant and $\delta_0$ is the Dirac mass at $0$, $h$ is a positive and smooth function. If $\alpha$ is not an integer the profile of bubbling solutions is so much easier but the corresponding equation does not reflect the most interesting applications the Liouville equations represent. If $\alpha$ is a positive integer, we say the singular source is quantized, the profile of bubbling solutions around a quantized singular source is the most challenging situation which also has the most significant applications. The goal in this article is to provide a complete understanding of this situation.

 It is known from the works of Kuo-Lin \cite{kuo-lin-jdg}, Bartolucci-Tarantello \cite{bart3} that if $\alpha\not \in \mathbb N$ ( the set of natural numbers) blowup solutions satisfy spherical Harnack inequality \footnote{ A sequence of bubbling solutions satisfying spherical Harnack inequality means the oscillation of solutions on each fixed radius around the singular point is uniformly bounded.} around the singular source and the asymptotic behavior is relatively easy to understand. However, when the strength of the singular source is a multiple of $4\pi$ ($\alpha\in \mathbb N$), the so called ``non-simple blowup" phenomenon does occur, which means the bubbling solutions may not satisfy spherical Harnack inequality and multiple local maximums near the singular source could appear. In this article we establish new estimates for non-simple bubbling solutions. Since the analysis is carried out in a neighborhood of the singular source, we just require the domain to be a small neighborhood of the origin, so our assumption of bubbling solutions is as follows: Let $\tilde u_k$ be a sequence of solutions of
\begin{equation}\label{t-u-k}
\Delta \tilde u_k(x)+\tilde h_k(x)e^{\tilde u_k}=4\pi N\delta_0, \quad \mbox{in}\quad B_{\tau}
\end{equation}
for some $\tau>0$ independent of $k$. $B_{\tau}$ is the ball centered at the origin with radius $\tau$.  In addition we postulate the usual assumptions on $\tilde u_k$ and $\tilde h_k$:
For a positive constant $C$ independent of $k$, the following holds:
\begin{equation}\label{assumption-1}
\left\{\begin{array}{ll}
\|\tilde h_k\|_{C^3(\bar B_{\tau})}\le C, \quad \frac 1C\le \tilde  h_k(x)\le C, \quad x\in \bar B_{\tau}, \\ \\
\int_{B_{\tau}} \tilde h_k e^{\tilde u_k}\le C,\\  \\
|\tilde u_k(x)- \tilde u_k(y)|\le C, \quad \forall x,y\in \partial B_{\tau},
\end{array}
\right.
\end{equation}
and since we study the asymptotic behavior of blowup solutions around the singular source, we assume that there is no blowup point except at the origin:
\begin{equation}\label{assump-2}
\max_{K\subset\subset B_{\tau}\setminus \{0\}} \tilde u_k\le C(K).
\end{equation}
Also,
for the convenience of notation we assume $\tilde h_k(0)=1$ and use the value of $\tilde u_k$ on $\partial B_{\tau}$ to define a harmonic function $\phi_k(x)$:
\begin{equation}\label{phi-k}
\left\{\begin{array}{ll}
\Delta \phi_k(x)=0,\quad \mbox{in}\quad B_{\tau},\\ \\
\phi_k(x)=u_k(x)-\frac{1}{2\pi \tau}\int_{\partial B_{\tau}}\tilde u_kdS,\quad x\in \partial B_{\tau}.
\end{array}
\right.
\end{equation}
Using the fact that $\Delta (\frac 1{2\pi}\log |x|)=\delta_0$, we set
\begin{equation}\label{uk-d}
u_k(x)=\tilde u_k(x)-2N\log |x|-\phi_k(x),
\end{equation}
which satisfies
\begin{equation}\label{eq-uk}
\Delta u_k(x)+|x|^{2N}h_k(x)e^{u_k}=0,\quad \mbox{ in }\quad B_{\tau}
\end{equation}
for
\begin{equation}\label{hk-d}
h_k(x)=\tilde h_k(x)e^{\phi_k(x)}.
\end{equation}
It is easy to see that $\phi_k(0)=0$ and $u_k$ is a constant on $\partial B_{\tau}$.

In this article we consider the case that:
\begin{equation}\label{no-sp-h}
\max_{x\in B_{\tau}} u_k(x)+2(1+N)\log |x|\to \infty,
\end{equation}
which is equivalent to saying that the spherical Harnack inequality does not hold for $u_k$.
It is well known \cite{kuo-lin-jdg, bart3} that $ u_k$ exhibits a non-simple blowup profile.  It is established in \cite{kuo-lin-jdg,bart3} that there are $N+1$ local maximum points of $ u_k$: $p_0^k$,....,$p_{N}^k$ and they are evenly distributed on $\mathbb S^1$ after scaling according to their magnitude: Suppose along a subsequence
$$\lim_{k\to \infty}p_0^k/|p_0^k|=e^{i\theta_0}, $$
then
$$\lim_{k\to \infty} \frac{p_l^k}{|p_0^k|}=e^{i(\theta_0+\frac{2\pi l}{N+1})}, \quad l=1,...,N. $$
For many reasons it is convenient to denote $|p_0^k|$ as $\delta_k$ and define $\mu_k$ as follows:
\begin{equation}\label{muk-dk}
\delta_k=|p_0^k|\quad \mbox{and }\quad \mu_k= u_k(p_0^k)+2(1+N)\log \delta_k.
\end{equation}

Since $p_l^k$'s are evenly distributed
around $\partial B_{\delta_k}$, standard results for Liouville equations around a regular blowup point can be applied to have $ u_k(p_l^k)= u_k(p_0^k)+o(1)$. Also, (\ref{no-sp-h}) gives $\mu_k\to \infty$. The interested readers may look into \cite{kuo-lin-jdg,bart3} for more detailed information.

In the first main result we prove a surprising vanishing type estimates for the first derivatives of the coefficient function $\log h_k$:
\begin{thm}\label{Van-n2} Let $u_k$, $\phi_k$, $h_k$, $\delta_k$, $\mu_k$ be defined by (\ref{eq-uk}), (\ref{phi-k}), (\ref{hk-d}), (\ref{muk-dk}) respectively. Then
\begin{equation}\label{first-order-v}
|\nabla \log h_k(0)|=O(\delta_k)+O(\delta_k^{-1}e^{-\mu_k}\mu_k).
\end{equation}
\end{thm}
Here we observe that if $\mu_ke^{-\mu_k}=o(\delta_k)$, the origin has to be a critical point of $h_k$. We can also prove a vanishing type theorem for second order derivatives of $h_k(0)$.
\begin{thm}\label{Van-n1}Let $u_k$, $\phi_k$, $h_k$, $\delta_k$, $\mu_k$ be defined by (\ref{eq-uk}), (\ref{phi-k}), (\ref{hk-d}), (\ref{muk-dk}) respectively. Then
$$\Delta (\log h_k)(0)=O(\delta_k)+O(\delta_k^{-2}\mu_ke^{-\mu_k}),\quad \mbox{ if }\quad N\ge 2, $$
and for $N=1$ we have
\begin{align}\label{thm-2-n-1}
&\partial_{e_ke_k}(\log h_k)(0)=O(\delta_k^{-2}\mu_ke^{-\mu_k})+O(\delta_k)\\
&\partial_{e_ke_k^{\perp}}(\log h_k)(0)=O(\delta_k^{-2}\mu_ke^{-\mu_k})+O(\delta_k), \nonumber
\end{align}
where $e_k=p_0^k/|p_0^k|$, $e_k^{\perp}$ is a unit vector orthogonal to $e_k$.
\end{thm}

As an important application of the vanishing type theorems we can address the question of approximating non-simple bubbling solutions by global solutions:
\begin{equation}\label{global-s}
\Delta U+|x|^{2N}e^{U}=0,\quad \mbox{in}\quad\mathbb R^2,\quad \int_{\mathbb R^2}|x|^{2N}e^{U}<\infty.
\end{equation}
 For regular Liouville equation, this type of approximation, initiated by Y.Y.Li \cite{licmp}, and further extended and refined by a series of works \cite{bclt,bart2,bart3,bart4,chenlin1,gluck,zhangcmp,zhangccm} played an important role in a number of applications such as degree counting theorems \cite{chenlin1,chenlin2}, uniqueness results \cite{bjly-1,bjly}, etc.  Our Theorem \ref{main-thm-2} below seems to be the first such result for quantized singular sources:

\begin{thm}\label{main-thm-2} Let $u_k$, $\phi_k$, $h_k$, $\delta_k$, $\mu_k$ be defined by (\ref{eq-uk}), (\ref{phi-k}), (\ref{hk-d}), (\ref{muk-dk}) respectively. If $\delta_k=O(e^{-\mu_k/4})$ we have, for some $c_1>0$ independent of $k$ and a sequence of global solutions $U_k$ of (\ref{global-s}) that
$$| u_k(x)-U_k(x)|\le c_1 \quad x\in B_{\tau}.$$
 \end{thm}
\medskip

\begin{rem} For $dist(x,0)\sim 1$, $u_k(x)=-u_k(p_0^k)+O(1)$. This is already established in \cite{kuo-lin-jdg, bart3}.
\end{rem}

\begin{rem} Theorem \ref{Van-n2} and Theorem \ref{Van-n1} are very surprising because it seems to be a common understanding among experts that it is not possible to obtain vanishing theorems around a singular blowup point. This is certainly the case when the singular source is not quantized \cite{zhangccm}. Theorems \ref{Van-n2} and \ref{Van-n1} not only provided vanishing estimates for first derivatives of coefficient function, but also those for the second derivatives. We would also like to point out that some vanishing estimates for bubbling solutions of Toda systems have been obtained exactly at singular sources \cite{lwz-aim}, \cite{zhang-imrn}.
\end{rem}

\begin{rem} Theorems \ref{Van-n2},\ref{Van-n1} and \ref{main-thm-2} also provide a dichotomy: If $\delta_k$ is too small to have vanishing rates in (\ref{first-order-v}) (\ref{thm-2-n-1}) useful, $v_k$ can be accurately approximated by global solutions. On the other hand, as long as $\delta_k$ is not extremely small, vanishing estimates imply that the blowup point has to be a critical point of the coefficient function.
\end{rem}

\begin{rem} Lemma 9 of \cite{kuo-lin-jdg} asserts that
$\delta_k^2=c\mu_ke^{-\mu_k}(1+o(1))$. However we found (4.18) of \cite{kuo-lin-jdg} incorrect. In fact there should not be any deterministic relation between $\mu_k$ and $\delta_k$, because for any $\xi_k\in \mathbb R^2$ and any $\lambda_k\in \mathbb R$,
$$U_k(x)=\log \frac{e^{\lambda_k}}{(1+\frac{e^{\lambda_k}}{8(1+N)^2}|x^{N+1}-\xi_k|^2)^2}$$
is a sequence of solutions to $\Delta U_k+|x|^{2N}e^{U_k}=0$.
\end{rem}

The non-simple bubbling situation and vanishing theorems proved in this article are closely related to so many important equations in geometry and physics that it is impossible to elaborate even a portion of them in reasonable detail. So we shall only try to explain a few representatives.
First for the following mean field equation defined on a Riemann surface $(M,g)$:
\begin{equation}\label{mean-sin}
\Delta_gu+\rho(\frac{h(x)e^{u(x)}}{\int_Mhe^{u}}-\frac 1{Vol_g(M)})=4\pi\sum_j \alpha_j (\delta_{p_j}-\frac 1{Vol_g(M)}),
\end{equation}
the solution $u$ represents a conformal metric with prescribed conic singularities (see \cite{erem-3,tr-1,tr-2}) . In particular if the singular source is quantized, the Liouville equation has close ties with Algebraic geometry, integrable system and number theory. The study of bubbling solutions near a quantized singular source reveals crucial information about
the curvature, topology of the manifold and related degree-counting formula. If we consider a slight variation of (\ref{mean-sin}) on a flat torus $T$,
\begin{equation}\label{flat-liou}
\Delta_g u+e^u=\rho \delta_0, \quad \mbox{ in }\quad T
\end{equation}
the solution $u$ of (\ref{flat-liou}) is related to the following complex valued Monge-Ampere equation
\begin{equation}\label{d-c-mp}
det(\frac{\partial^2w}{\partial z_i\partial z_j})^d=e^{-w}, \quad \mbox{on}\quad (T\setminus \{0\})^d,
\end{equation}
as follows: Let $w(z_1,...,z_d)=-\sum_{i=1}^du(z_i)+d\log 4$, then $w$ satisfies (\ref{d-c-mp}) with a logarithmic singularity along the normal crossing divisor
$D=T^d\setminus (T\setminus \{0\})^d$. Also, the bubbling solutions of (\ref{flat-liou}) can be translated as blowup solutions of (\ref{d-c-mp}). The main theorems of this article could shed light to the understanding of the geometry related the degenerate complex Monge-Amepere equation like (\ref{d-c-mp}).

In physics the mean field equation (\ref{mean-sin})
 is derived from the mean field limit of point vortices in the Euler flow \cite{caglioti-1,caglioti-2} and serves as a model equation in the Chern-Simons-Higgs theory \cite{jackiw} and in the electroweak theory \cite{ambjorn}, etc. Among many equations related to the mean field equation we mention one in particular:
 the self-dual Bogomolny type equation:
 \begin{equation}\label{phy-open}
 \left\{\begin{array}{ll}
 \Delta u+4g^2e^u+g^2e^w=4\pi\sum_{l=1}^mn_l\delta_{p_l},\quad \mbox{ in } \quad T\\
 \Delta w-2g^2e^u-\frac{g^2}{2\cos^2\theta}(e^w-\phi_0^2)=0,
 \end{array}
 \right.
 \end{equation}
 where $\phi_0$,$\theta$, $g$ are constants. Integrating the equation we have
 $$4g^2\int_Te^u=\frac{4\pi N-g^2\phi_0^2|T|}{\sin^2\theta}$$
 and
 $$g^2\int_Te^{w}=\frac{g^2\phi_0^2|T|-4\pi \cos^2\theta N}{\sin^2\theta} $$
 where $N=\sum_{l=1}^mn_l$. In order for (\ref{phy-open}) to have a solution, a necessary condition is
 \begin{equation}\label{phy-sur}
 g^2\phi_0^2<\frac{4\pi T}{|T|}<\frac{g^2\phi_0^2}{\cos^2\theta}
 \end{equation}
 Yang \cite{yang} posted an important question of whether (\ref{phy-sur}) is also sufficient for the existence of a solution to (\ref{phy-open}). Such a conjecture was partially solved by Bartolucci-De Marchis \cite{bart-de} and Chen-Lin \cite{chen-lin-last-cpam} if
$\frac{4\pi N-g^2\phi_0^2|T|}{\sin^2\theta}\not \in 8\pi \mathbb N$. Thus the most interesting open case of the conjecture is when
$\frac{4\pi N-g^2\phi_0^2|T|}{\sin^2\theta}$ is a multiple of $8\pi$. The main results of this article could be crucial for the complete resolution of this conjecture.

The phenomena of non-simple bubbling solutions not only occur in single equations, but also in systems. In a recent work of the second author and Gu \cite{gu-zhang}, the non-simple blowup solutions are studied for singular Liouville systems.

To end the introduction we would like to briefly explain the idea of the proof. In \cite{kuo-lin-jdg} and \cite{bart3} it is already established that there are exactly $N+1$ local maximum points evenly distributed around the origin. Kuo-Lin and Bartolucci-Tarantello independently obtained this important information by studying the Pohozaev identity around each local maximum. In this  article our investigation starts from Kuo-Lin and Bartolucci-Tarantello's result and go way deeper in our analysis. First we prove that all the local maximum points of bubbling solutions are perturbations of roots of a unity. All these perturbations are determined by the logarithmic of a coefficient function. In order to fully understand the nature of the problem, we launch a thorough investigation of the algebraic structure of the bubbling solutions and prove that the first derivatives of the coefficient function must vanish. Not only so the second derivatives also vanish at certain speed at the singular source. Since the singular Liouville equation with quantized singularity appears in many context, we expect both the result and the approach of this article to have impact to related problems.

The organization of this article is as follows: In section two we establish some preliminary estimates, in section three we establish rough locations of local maximum points. In this section the main matrix that contains most algebraic structure of this article is determined. Then in section four we prove the vanishing estimate for first order derivatives of the coefficient function.  The approximation by global solutions (Theorem \ref{main-thm-2}) is proved in section five and the proof of vanishing theorem of second derivatives (Theorem \ref{Van-n1}) is arranged in
section six. In appendix A we prove a sharp estimate of bubbling solutions if the spherical Harnack inequality holds, in this section we employ the idea in \cite{lwz-jems} to defeat the lack of symmetry. Finally in Appendix B we cite Gluck's result \cite{gluck} for a sample computation that is used multiple times in this article.

\medskip
{\bf Notation:} We will use $B(x_0,r)$ to denote a ball centered at $x_0$ with radius $r$. If $x_0$ is the origin we use $B_r$. $C$ represents a positive constant that may change from place to place.

\medskip
{\bf Acknowledgement:} We are grateful to Wen Yang for stimulating discussions and anonymous referees for many insightful remarks and excellent suggestions.

\section{Preliminary Discussions}

First we recall that $|p_0^k|=\delta_k$, so we write $p_0^k$ as $p_0^k=\delta_ke^{i\theta_k}$ and define $v_k$ as
\begin{equation}\label{v-k-d}
v_k(y)=u_k(\delta_k ye^{i\theta_k})+2(N+1)\log \delta_k,\quad |y|<\tau \delta_k^{-1}.
\end{equation}
If we write out each component, (\ref{v-k-d}) is
$$
v_k(y_1,y_2)=u_k(\delta_k(y_1\cos\theta_k-y_2\sin\theta_k),\delta_k(y_1\sin\theta_k+y_2\cos\theta_k))+2(1+N)\log \delta_k. $$
Then it is standard to verify that $v_k$ solves

\begin{equation}\label{e-f-vk}
\Delta v_k(y)+|y|^{2N}\mathfrak{h}_k(\delta_k y)e^{v_k(y)}=0,\quad |y|<\tau/\delta_k,
\end{equation}
where
\begin{equation}\label{frak-h}
\mathfrak{h}_k(x)=h_k(xe^{i\theta_k}),\quad |x|<\tau.
\end{equation}
Thus the image of $p_0^k$ after scaling is $Q_1^k=e_1=(1,0)$.
Let $Q_1^k$, $Q_2^k$,...,$Q_{N}^k$ be the images of $p_i^k$ $(i=1,...,N)$ after the scaling:
$$Q_l^k=\frac{p_l^k}{\delta_k}e^{-i\theta_k},\quad l=1,...,N. $$
 It is established by Kuo-Lin in \cite{kuo-lin-jdg} and independently by Bartolucci-Tarantello in \cite{bart3} that
\begin{equation}\label{limit-q}
\lim_{k\to \infty} Q_l^k=\lim_{k\to \infty}p_l^k/\delta_k=e^{\frac{2l\pi i}{N+1}},\quad l=0,....,N.
\end{equation}

Choosing $\epsilon>0$ small and independent of $k$, we can make disks centered at $Q_l^k$ with radius $\epsilon$ (denoted as $B(Q_l^k,\epsilon ) $) mutually disjoint. Let
\begin{equation}\label{v-muk}
\mu_k=\max_{B(Q_0^k,\epsilon)} v_k.
\end{equation}
Since $Q_l^k$ are evenly distributed around $\partial B_1$, it is easy to use standard estimates for single Liouville equations (\cite{zhangcmp,gluck,chenlin1}) to obtain
$$\max_{B(Q_l^k,\epsilon)}v_k=\mu_k+o(1),\quad l=1,...,N. $$

Recall that $v_k$ satisfies (\ref{e-f-vk})
and $v_k$ is a constant on $\partial B(0,\tau \delta_k^{-1})$.  The Green's representation formula for $v_k$ gives,
$$v_k(y)=\int_{\Omega_k}G_k(y,\eta) |\eta |^{2N} \mathfrak{h}_k(\delta_k\eta)e^{v_k(\eta)}d\eta +v_k|_{\partial \Omega_k}
$$
where $\Omega_k=B(0,\tau \delta_k^{-1})$ and
\begin{equation}\label{green-e}
G_k(y,\eta)=-\frac{1}{2\pi}\log |y-\eta |+H_k(y,\eta)
\end{equation}
where
$$H_k(y,\eta)=\frac{1}{2\pi}\log \bigg (\frac{|\eta |}{\tau \delta_k^{-1}}|\frac{\tau^2\delta_k^{-2}\eta}{|\eta |^2}-y|\bigg ). $$

Also for $r>2$, let $\bar v_k(r)$ be the spherical average of $v_k$ on $\partial B_r$, then we have
$$\frac{d}{dr}\bar v_k(r)=\frac{d}{dr}\bigg (\frac{1}{2\pi r}\int_{B_r}\Delta v_k\bigg )=-\frac{8(N+1)\pi+o(1)}{2\pi r}. $$
Because of the fast decay of $\bar v_k(r)$ it is easy to use the Green's representation of $v_k$ to obtain the following stronger estimate of $v_k$:
\begin{equation}\label{vk-crude}
v_k(y)=- \mu_k-(4N+4)\log |y|+O(1),\quad 2<|y|<\tau \delta_k^{-1}.
\end{equation}

Now we consider $v_k$ around $Q_l^k$. Using the results in \cite{chenlin1,zhangcmp,gluck} we have, for $v_k$ in $B(Q_l^k,\epsilon )$, the following gradient estimate:

 \begin{equation}\label{gra-each-p}
\delta_k \nabla  (\log \mathfrak{h}_k)(\delta_k \tilde Q_l^k)+2N\frac{\tilde Q_l^k}{|\tilde Q_l^k|^2}+\nabla \phi_l^k(\tilde Q_l^k)=O( \mu_k e^{- \mu_k}),
\end{equation}
where $\phi_l^k$ is the harmonic function that eliminates the oscillation of $v_k$ on $\partial B(Q_l^k,\epsilon)$ and $\tilde Q_l^k$ is the maximum of $v_k-\phi_l^k$ that satisfies
\begin{equation}\label{close-2}
\tilde Q_l^k-Q_l^k=O( e^{-\mu_k}).
\end{equation}
Using (\ref{close-2}) in (\ref{gra-each-p}) we have
\begin{equation}\label{pi-each-p}
\delta_k \nabla  (\log \mathfrak{h}_k)(\delta_k  Q_l^k)+2N\frac{Q_l^k}{| Q_l^k|^2}+\nabla \phi_l^k( Q_l^k)=O( \mu_k e^{- \mu_k}).
\end{equation}
For the discussion in this section we use the following version of (\ref{pi-each-p}):
\begin{equation}\label{location-delta}
\delta_k\nabla  (\log \mathfrak{h}_k)(0)+2N\frac{Q_l^k}{|Q_l^k|^2}+\nabla \phi_l^k(Q_l^k)=O(\delta_k^2)+O( \mu_k e^{- \mu_k})
\end{equation}
and the first estimate of $\nabla \phi_l^k(Q_l^k)$ is
\begin{lem}\label{phi-k-e} For $l=0,...,N$,
\begin{equation}\label{gra-each-2}
\nabla \phi_l^k(Q_l^k)
=-4\sum_{m=0,m\neq l}^N\frac{Q_l^k-Q_m^k}{|Q_l^k-Q_m^k|^2}
+E
\end{equation}
where
\begin{equation}\label{s-err}
E=O(\delta_k^2)+O(\mu_k e^{-\mu_k}).
\end{equation}
\end{lem}

\noindent{\bf Proof of Lemma \ref{phi-k-e}:}

\medskip

From the expression of $v_k$ on $\Omega_k=B(0,\tau \delta_k^{-1})$ we have, for $y$ away from bubbling disks,
\begin{align}\label{pi-e-1}
v_k(y)&=v_k|_{\partial \Omega_k}+\int_{\Omega_k}G(y,\eta)|\eta |^{2N}\mathfrak{h}_k(\delta_k\eta)e^{v_k(\eta)}d\eta\\
&=v_k|_{\partial \Omega_k}+\sum_{l=0}^{N}G(y,Q_l^k)\int_{B(Q_l,\epsilon)}|\eta |^{2N}\mathfrak{h}_k(\delta_k\eta)e^{v_k}d\eta \nonumber\\
&+\sum_l\int_{B(Q_l,\epsilon)}(G(y,\eta)-G(y,Q_l^k))|\eta |^{2N}\mathfrak{h}_k(\delta_k\eta)e^{v_k}d\eta+O(\mu_ke^{-\mu_k}). \nonumber
\end{align}
Before we evaluate each term, we use a sample computation which will be used repeatedly: Suppose $f$ is a smooth function defined on $B(Q_0^k,\epsilon)$, then
we evaluate
\begin{equation}\label{sample-c}
\int_{B(Q_0^k,\epsilon)}f(\eta)|\eta |^{2N}\mathfrak{h}_k(\delta_k \eta)e^{v_k(\eta)}d\eta.
\end{equation}
Let $\tilde Q_0^k$ be the maximum of $v_k-\phi_0^k$, and set $\hat h_k(y)=|y|^{2N}\mathfrak{h}_k(\delta_ky)$ then it is known \cite{zhangcmp,gluck} that
\begin{equation}\label{close-Q}
\tilde Q_0^k-Q_0^k=O(e^{-\mu_k}).
\end{equation}
Moreover, it is derived that
\begin{align}\label{gluck-t}
v_k(y)=\phi_0^k(y)+\log \frac{e^{\mu_k}}{(1+e^{\mu_k}\frac{\hat h_k(\tilde Q_0^k)}{8}|y-\tilde Q_0^k|^2)^2}\\
+c_ke^{-\mu_k}(\log (2+e^{u_k/2}|y-\tilde Q_0^k|))^2+O(\mu_ke^{-\mu_1^k}). \nonumber
\end{align}
We use $U_k$ to denote the leading global solution.
Using the expansion above and symmetry we have
\begin{align}\label{sampled-cm}
&\int_{B(\tilde Q_0^k,\epsilon)}f(\eta)|\eta |^{2N}\mathfrak{h}_k(\delta_k \eta)e^{v_k(\eta)}d\eta\\
=&\int_{B(\tilde Q_0^k,\epsilon)}f(\eta)\hat h_k(y)e^{v_k(y)}dy\nonumber \\
=&\int(f(\tilde Q_0^k)+\nabla f(\tilde Q_0^k)\tilde \eta+O(|\tilde \eta|^2)(\hat h_k(\tilde Q_0^k)+\nabla \hat h_k(\tilde Q_0^k)(\tilde \eta)
+O(|\tilde \eta|^2)\nonumber \\
&\cdot e^{U_k}(1+\phi_0^k+c_ke^{-\mu_k}(\log (2+e^{\mu_k/2}|\tilde \eta|))^2+O(\mu_ke^{-\mu_k})+O(|\tilde \eta|^2))d\eta \nonumber
\end{align}
where $\tilde \eta=\eta-\tilde Q_0^k$.  In the evaluation we use symmetry, for example,
$$\int_{B(\tilde Q_0^k,\epsilon)}e^{U_k}\tilde \eta d\eta=0. $$ Also
using $\phi_0^k(Q_0^k)=0$ and $Q_0^k-\tilde Q_0^k=O(e^{-\mu_k})$ we have
$$\int_{B(\tilde Q_0^k,\epsilon)} e^{U_k}\phi_k=O(e^{-\mu_k}). $$
Carrying out computations in (\ref{sampled-cm}) we arrive at
$$\int_{B(Q_0^k,\epsilon)}f(\eta)|\eta |^{2N}\mathfrak{h}_k(\delta_k \eta)e^{v_k(\eta)}d\eta=8\pi f(\tilde Q_0^k)+O(\mu_k e^{-\mu_k}). $$
Since (\ref{close-Q}) holds we further have
\begin{equation}\label{sample-2}
\int_{B(Q_0^k,\epsilon)}f(\eta)|\eta |^{2N}\mathfrak{h}_k(\delta_k \eta)e^{v_k(\eta)}d\eta=8\pi f( Q_0^k)+O(\mu_k e^{-\mu_k}).
\end{equation}

Using the method of (\ref{sample-2}) in the evaluation of each term in (\ref{pi-e-1}) we have,
$$
v_k(y)
=v_k|_{\partial \Omega_k}-4\sum_{l=0}^{N}\log |y- Q_l^k|
+8\pi\sum_{l=0}^{N}H(y,Q_l^k)+O(\mu_ke^{-\mu_k}). $$

The harmonic function that kills the oscillation of $v_k$ around $Q_m^k$ is
\begin{align*}
\phi_m^k=-4\sum_{l=0,l\neq m}^N(\log |y-Q_l^k|-\log |Q_m^k-Q_l^k|)\\
+8\pi\sum_{l=0}^{N}(H(y,Q_l^k)-H(Q_m^k,Q_l^k))+O(\mu_ke^{-\mu_k}).
\end{align*}

The corresponding estimate for $\nabla \phi_m^k$ is
$$
\nabla \phi_m^k(Q_m^k)=-4\sum_{l=0,l\neq m}^N\frac{Q_m^k-Q_l^k}{|Q_m^k-Q_l^k|^2}
+8\pi\sum_{l=0}^{N}\nabla_1 H(Q_m^k,Q_l^k)+O(\mu_ke^{-\mu_k}).
$$
where $\nabla_1$ stands for the differentiation with respect to the first component. From the expression of $H$, we have
\begin{align}\label{grad-H}
\nabla_1H(Q_m^k,Q_l^k)&=\frac 1{2\pi}\frac{Q_m^k-\tau^2\delta_k^{-2}Q_l^k/|Q_l^k|^2}{|Q_m^k-\tau^2\delta_k^{-2}Q_l^k/|Q_l^k|^2|^2}\\
&=\frac 1{2\pi}\tau^{-2}\delta_k^2\frac{\tau^{-2}\delta_k^2Q_m^k-Q_l^k/|Q_l^k|^2}{|Q_l^k/|Q_l^k|^2-\tau^{-2}\delta_k^2Q_m^k|^2}\nonumber \\
&=-\frac{1}{2\pi}\tau^{-2}\delta_k^2e^{\frac{2\pi i l}{N+1}}+O(\sigma_k\delta_k^2). \nonumber
\end{align}
where $\sigma_k=\max_l|Q_l^k-e^{\frac{2\pi i l}{N+1}}|$. Later we shall obtain more specific estimate of $\sigma_k$.
Thus
\begin{align}\label{phimk-g}
&\nabla \phi_m^k(Q_m^k) \\
=&-4\sum_{l=0,l\neq m}^N\frac{Q_m^k-Q_l^k}{|Q_m^k-Q_l^k|^2}
-4\tau^{-2}\delta_k^2\sum_{l=0}^{N}e^{\frac{2\pi il}{N+1}}+O(\sigma_k\delta_k^2)+O(\mu_ke^{-\mu_k}) \nonumber \\
=&-4\sum_{l=0,l\neq m}^N\frac{Q_m^k-Q_l^k}{|Q_m^k-Q_l^k|^2}
+O(\sigma_k\delta_k^2)+O(\mu_ke^{-\mu_k}) \nonumber
\end{align}
where we have used $\sum_{l=0}^{N}e^{2\pi l i/(N+1)}=0$.
Since we don't have the estimate of $\sigma_k$ now we have
$$
\nabla \phi_m^k(Q_m^k)=-4\sum_{l=0,l\neq m}^N\frac{Q_m^k-Q_l^k}{|Q_m^k-Q_l^k|^2}+E
$$
Lemma \ref{phi-k-e} is established. $\Box$

\section{Location of blowup points}

In this section we establish a first description of the locations of $Q_l^k$ ($l=0,..,N)$). The result of this section will be used later to obtain vanishing estimates of the coefficient function $h_k$.

Let $E$ be defined as in (\ref{s-err}).  The Pohozaev identity around $Q_l^k$ now reads

$$-4\sum_{j=0,j\neq l}^N\frac{Q_l^k-Q_j^k}{|Q_l^k-Q_j^k|^2}+2N\frac{Q_l^k}{|Q_l^k|^2}=-\nabla (\log \mathfrak{h}_k)(0)\delta_k +E. $$
Using $L_k$ to denote $\nabla (\log  \mathfrak{h}_k)(0)$, we have, treating every term as a complex number,
$$N\frac{1}{\bar Q_l^k}=2\sum_{j=0,j\neq l}^N\frac{1}{\bar Q_l^k-\bar Q_j^k}-\frac{L_k}2\delta_k+E, $$
where $\bar Q_l^k$ is the conjugate of $Q_l^k$.
Thus
\begin{equation}\label{pi-N}
N=2\sum_{j=0,j\neq l}^N\frac{Q_l^k}{Q_l^k-Q_j^k}-\frac{\bar L_k}2\delta_k Q_l^k+E.
\end{equation}
Let $\beta_l=2\pi l/(N+1)$, we write $Q_l^k=e^{i\beta_l }+q_l^k$ for $q_l^k\to 0$. Then we write the first term on the right hand side of (\ref{pi-N}) as
\begin{align*}
&\frac{Q_l^k}{Q_l^k-Q_j^k}=\frac{e^{i\beta_l}+q_l^k}{e^{i\beta_l}-e^{i\beta_j}+q_l^k-q_j^k}\\
=&\frac{e^{i\beta_l}+q_l^k}{(e^{i\beta_l}-e^{i\beta_j})(1+(q_l^k-q_j^k)/(e^{i\beta_l}-e^{i\beta_j}))}\\
=&\frac{e^{i\beta_l}}{e^{i\beta_l}-e^{i\beta_j}}+\frac{q_l^k}{e^{i\beta_l}-e^{i\beta_j}}-
\frac{e^{i\beta_l}}{(e^{i\beta_l}-e^{i\beta_j})^2}(q_l^k-q_j^k)+O(\sigma_k^2)\\
=&\frac{e^{i\beta_l}}{e^{i\beta_l}-e^{i\beta_j}}+\frac{e^{i\beta_l}q_j^k-e^{i\beta_j}q_l^k}{(e^{i\beta_l}-e^{i\beta_j})^2}+O(\sigma_k^2).
\end{align*}

Using
\begin{equation}\label{e-N}
N=2\sum_{j=0,j\neq l}^N\frac{e^{i\beta_l}}{e^{i\beta_l}-e^{i\beta_j}},
\end{equation}
we write (\ref{pi-N}) as
\begin{equation}\label{var-1}
\sum_{j=0,j\neq l}^N\frac{e^{i\beta_l}q_j^k-e^{i\beta_j}q_l^k}{(e^{i\beta_l}-e^{i\beta_j})^2}-\frac{\bar L_k}{4}\delta_k e^{i\beta_l}=E+O(\sigma_k^2)
\end{equation}
for $l=0,1,2,....,N$.
For convenience we set
$$q_l^k=e^{i\beta_l}m_l^k \quad \mbox{ and  }\quad \beta_{jl}=\beta_j-\beta_l $$
to  reduce (\ref{var-1}) to
\begin{align}\label{main-1}
&\sum_{j=0,j\neq l}^N\frac{e^{i\beta_{jl}}m_j^k}{(1-e^{i\beta_{jl}})^2}-\bigg (\sum_{j=0,j\neq l}^N\frac{e^{i\beta_{jl}}}{(1-e^{i\beta_{jl}})^2}\bigg )m_l^k
-\frac{\bar L_k}4\delta_k e^{i\beta_l}\\
&=E+O(\sigma_k^2)+O(\delta_k \sigma_k) \nonumber
\end{align}
for $l=0,1.....,N$.
It is easy to verify that
\begin{equation}\label{trig-i-1}
\frac{e^{i\theta}}{(1-e^{i\theta})^2}=\frac{1}{2(\cos \theta -1)}=(-\frac 14)\frac{1}{\sin^2(\theta/2)}.
\end{equation}
To deal with coefficients of $m_j^k$ in (\ref{main-1}) we set
$$d_j=\frac{1}{\sin^2(\frac{j\pi}{N+1})},\quad j=1,...,N $$
and
$$D=\sum_{j=0,j\neq l}^N d_{|j-l|}. $$
Since $d_l=d_{N+1-l}$ it is easy to check that $D$ does not depend on $l$:
\begin{equation}\label{D-guess}
D=\sum_{k=1}^N d_k=\sum_{k=1}^N\frac{1}{\sin^2(\frac{k\pi}{N+1})}=\frac{N^2+2N}{3}.
\end{equation}
Now (\ref{main-1}) can be written as
\begin{equation}\label{main-2}
-\sum_{j\neq l,j=0}^N d_{|j-l|}m_j^k+Dm_l^k-\bar{L_k}\delta_k e^{i\beta_l}=E+O(\sigma_k^2),\quad l=0,....,N.
\end{equation}

For $l=0$, we have $\beta_0=0$ and $m_0^k=0$. Thus from (\ref{main-2}) we have
\begin{equation}\label{main-3}
-\sum_{j=1}^N d_jm_j^k-\bar L_k\delta_k =E+O(\sigma_k^2),
\end{equation}
If we take $(m^k_1,...,m^k_n)$ as unknowns in (\ref{main-2}), the last $N$ equations of (\ref{main-2}) ( for $l=1,...,N$) can be written as
\begin{equation}\label{e-m}
A\left(\begin{array}{c}
m^k_1\\
m^k_2\\
\vdots\\
m^k_N
\end{array}
\right)=\bar L_k \delta_k \left(\begin{array}{c}
e^{i\beta_1}\\
e^{i\beta_2}\\
\vdots\\
e^{i\beta_N}
\end{array}
\right)+E+O(\sigma_k^2).
\end{equation}
where
$$A=\left(\begin{array}{cccc}
D & -d_1 & ... & -d_{N-1} \\
-d_1 & D & ... & -d_{N-2} \\
\vdots & \vdots & ... & \vdots \\
-d_{N-1} & -d_{N-2} & ... & D
\end{array}
\right ) $$

Since $D=|d_1|+...+|d_N|$ and each $d_i>0$, we see that the matrix is invertible. Note that the magnitude of $|m_i^k|$ is the same as $\sigma_k$, thus
from equation (\ref{e-m}) we obtain
\begin{equation}\label{m-position}
\left(\begin{array}{c}
m_1^k\\
m_2^k\\
\vdots\\
m_N^k
\end{array}
\right)=A^{-1}\delta_k \bar L_K\left(\begin{array}{c}
e^{i\beta_1}\\
e^{i\beta_2}\\
\vdots\\
e^{i\beta_N}
\end{array}
\right)+E.
\end{equation}
With this fact we can further write $\nabla_1 H(Q_m^k,Q_l^k)$  in (\ref{grad-H}) as
\begin{equation}\label{grad-H-1}
\nabla_1H(Q_m^k,Q_l^k)
=-\frac{1}{2\pi}\tau^{-2}\delta_k^2e^{\frac{2\pi i l}{N+1}}+O(\delta_k^3)+O(\delta_k^2\mu_ke^{-\mu_k}),
\end{equation}
and $\nabla \phi_l^k(Q_l^k)$ in (\ref{gra-each-2}) and (\ref{phimk-g}) as
\begin{equation}\label{gra-each-3}
\nabla \phi_l^k(Q_l^k)
=-4\sum_{m\neq l}\frac{Q_l^k-Q_m^k}{|Q_l^k-Q_m^k|^2}
+O(\delta_k^3)+O(\mu_ke^{-\mu_k}).
\end{equation}

Using $(a^{ij})_{n\times n}$ to denote $A^{-1}$, we rewrite (\ref{m-position}) as
\begin{equation}\label{loca-m}
m_l^k=\delta_k\bar L_k\sum_{s=1}^na^{ls}e^{i\beta_s}+E,\quad l=1,...,N.
\end{equation}
This expression will be used in more refined analysis later.

%\begin{rem} In \cite{kuo-lin-jdg} the authors claimed that for non-simple blowup solutions, $\delta_k^2=c\mu_ke^{-\mu_k}(1+o(1))$. However we found that (4.18) of \cite{kuo-lin-jdg} is incorrect, which leads to a mistake in (4.26) and the estimate immediately after (4.26). As a result, the conclusion of Lemma 9 of \cite{kuo-lin-jdg}, as well as (1.26),(1.27), Corollary 1 and Theorem 5 are incorrect as well. But we firmly believe that except for what we just mentioned, the proof of all the other main results in \cite{kuo-lin-jdg} is ironclad.
%\end{rem}

\section{Vanishing estimates of $\nabla \mathfrak{h}_k(0)$}

In this section we establish a crucial estimate on the vanishing order of $\nabla (\log \mathfrak{h}_k)(0)$.   The main result in this section is Theorem \ref{Van-n2}.

\medskip

\noindent{\bf Proof of Theorem \ref{Van-n2}:} Before we get into the details we would like to mention that the most important observation is that the difference between the Pohozaev identities of $v_k$ and that of a global solution can be evaluated in detail. Using the result in the previous section we shall see that the coefficient of the desired first derivatives of $\mathfrak{h}_k(0)$ is not zero.

First we recall the equation for $v_k$:
$$\Delta v_k+\mathfrak{h}_k(\delta_ky)|y|^{2N}e^{v_k}=0,\quad |y|<\tau \delta_k^{-1} $$
with $v_k=$ constant on $\partial B(0,\tau \delta_k^{-1})$. Moreover $v_k(e_1)=\mu_k$. Now we set
$V_k$ to be the solution to
$$\Delta V_k+\mathfrak{h}_k(\delta_ke_1)|y|^{2N}e^{V_k}=0,\quad \mbox{in}\quad \mathbb R^2, \quad \int_{\mathbb R^2}|y|^{2N}e^{V_k}<\infty $$
such that $V_k$ has its local maximums at $e^{\frac{2\pi il}{N+1}}$ for $l=0,...,N$ and $V_k(e_1)=\mu_k$.
By the classification theorem of Prajapat-Tarantello \cite{prajapat}, the expression of $V_k$ is
\begin{equation}\label{Vk-exp}
V_k(y)=\log \frac{e^{\mu_k}}{(1+\frac{e^{\mu_k}\mathfrak{h}_k(\delta_ke_1)}{8(N+1)^2}|y^{N+1}-e_1|^2)^2}.
\end{equation}
In particular for $|y|\sim \delta_k^{-1}$,
$$V_k(y)=-\mu_k-4(N+1)\log \delta_k^{-1}+C+O(\delta_k^{N+1})+O(e^{-\mu_k}). $$

Let $w_k=v_k-V_k$ and $\Omega_k=B(0,\tau \delta_k^{-1})$,  we shall derive a precise, point-wise estimate of $w_k$ in $B_3\setminus \cup_{l=1}^{N}B(Q_l^k,\lambda)$ where $\lambda>0$ is a small number independent of $k$. Here we note that among $N+1$ local maximum points, we already have $e_1$ as a common local maximum point for both $v_k$ and $V_k$ and we shall prove that $w_k$ is very small in $B_3$ if we exclude all bubbling disks except the one around $e_1$. Before we carry out more specific computation we emphasize the importance of
\begin{equation}\label{control-e}
w_k(e_1)=|\nabla w_k(e_1)|=0.
\end{equation}
First we consider the Green's representation of $v_k$ on $\Omega_k$
$$v_k(y)=\int_{\Omega_k}G_k(y,\eta)|\eta|^{2N}\mathfrak{h}_k(\delta_k\eta)e^{v_k(\eta)}d\eta+v_k|_{\partial \Omega_k}, $$
where in $\Omega_k$, the Green's function $G_k(y,\eta)$ is written as in (\ref{green-e}).
Similarly for $V_k$ we obtain from the asymptotic expansion of $V_k$ to have
$$V_k(y)=\int_{\Omega_k}G_k(y,\eta)|\eta|^{2N}\mathfrak{h}_k(\delta_ke_1)e^{V_k(\eta)}d\eta+\bar V_k|_{\partial \Omega_k}+E, $$
where $\bar V_k|_{\partial \Omega_k}$ is the average of $V_k$ on $\partial \Omega_k$ and we have used the fact that the oscillation of $V_k$ on $\partial \Omega_k$ is $O(\delta_k^{N+1})+O(e^{-\mu_k})$.
The combination of these two equations gives, for $y\in B_3$,
$$w_k(y)=\int_{\Omega_k}G_k(y,\eta)|\eta|^{2N}(\mathfrak{h}_k(\delta_k\eta)e^{v_k}-\mathfrak{h}_k(\delta_ke_1)e^{V_k})d\eta+c_k+E. $$
Using $w_k(e_1)=0$, we clearly have
\begin{equation}\label{wk-pre-1}
w_k(y)=\int_{\Omega_k}(G_k(y,\eta)-G_k(e_1,\eta))|\eta |^{2N}(\mathfrak{h}_k(\delta_k\eta)e^{v_k}-\mathfrak{h}_k(\delta_ke_1)e^{V_k})d\eta+E.
\end{equation}
If we concentrate on $y\in B_3\setminus \cup_{l=0}^{N}B(Q_l^k,\epsilon)$ for $\epsilon>0$ small, we first claim that the harmonic part of $G_k$ only contributes $O(\delta_k^2)$ in the estimate:
\begin{equation}\label{har-minor}
\int_{\Omega_k}(H_k(y,\eta)-H_k(e_1,\eta))|\eta |^{2N}(\mathfrak{h}_k(\delta_k\eta)e^{v_k}-\mathfrak{h}_k(\delta_ke_1)e^{V_k})d\eta=O(\delta_k^2)
\end{equation}
where $G_k(y,\eta)=-\frac 1{2\pi}\log |y-\eta |+H_k(y,\eta)$. Indeed,
$$H_k(y,\eta)=\frac 1{2\pi}\log R_k+\frac{1}{2\pi}\log |\frac{y}{|y|}-\frac{\eta |y|}{R_k^2}|,\quad R_k=\tau \delta_k^{-1}. $$
Using this expression in the evaluation of the left hand side of (\ref{har-minor}), taking into consideration of the fast decay of $e^{V_k}$ and $e^{v_k}$, we obtain (\ref{har-minor}) with elementary estimates.

Combining (\ref{wk-pre-1}) and (\ref{har-minor}), we have
\begin{equation}\label{crucial-1}
w_k(y)=\int_{\Omega_k}\frac 1{2\pi}(\log |e_1-\eta|-\log |y-\eta|)|\eta |^{2N}(\mathfrak{h}_k(\delta_k\eta)e^{v_k}-\mathfrak{h}_k(\delta_ke_1)e^{V_k})d\eta+E.
\end{equation}

To evaluate the right hand side, first we observe that the integration away from the $N+1$ bubbling disks is $O(e^{-\mu_k}\delta_k)$ based on the asymptotic behavior of $V_k$. So we focus on the integration over $B(Q_l^k,\epsilon)$ for $l=0,..,N$.

If $l=0$, $Q_0^k=e_1$ by definition, we show that the integration over $B(e_1,\epsilon)$ is a very small error.

\medskip

For $l=1,...,N$, recall that the local maximum points of $V_k$ are at $e^{i\beta_l}$. By standard estimates for single Liouville equations, we have
\begin{align}\label{other-1}
&\frac 1{2\pi}\int_{B(Q_l^k,\epsilon)}(\log \frac{|e_1-\eta|}{|y-\eta|}|\eta|^{2N}(\mathfrak{h}_k(\delta_k\eta)e^{v_k}-\mathfrak{h}_k(\delta_ke_1)e^{V_k})d\eta \nonumber\\
=&4\log \frac{|e_1-Q_l^k|}{|y-Q_l^k|}-4\log \frac{|e_1-e^{i\beta_l}|}{|y-e^{i\beta_l}|}+O(\mu_ke^{-\mu_k})\nonumber\\
=&4\log \frac{|e_1-Q_l^k|}{|e_1-e^{i\beta_l}|}-4\log \frac{|y-Q_l^k|}{|y-e^{i\beta_l}|}+O(\mu_ke^{-\mu_k}).
\end{align}

Finally for $B(e_1,\epsilon)$ we find that by symmetry and the expansion of bubbles we have
\begin{align*}
\frac 1{2\pi}\int_{B(e_1,\lambda)}\log |e_1-\eta||\eta|^{2N}(\mathfrak{h}_k(\delta_k\eta)e^{v_k}-\mathfrak{h}_k(\delta_ke_1)e^{V_k})d\eta\\
=O(\delta_k\mu_k^2e^{-\mu_k})+E,
\end{align*}
and for $y\in B_3$ away from bubbling disks,
$$
\frac 1{2\pi}\int_{B(e_1,\lambda)}\log |y-\eta||\eta|^{2N}(\mathfrak{h}_k(\delta_k\eta)e^{v_k}-\mathfrak{h}_k(\delta_ke_1)e^{V_k})d\eta
=E.
$$

Putting the estimates in different regions together we have
\begin{equation}\label{e-wk-important}
w_k(y)=\sum_{l=1}^{N}(4\log \frac{|e_1-Q_l^k|}{|e_1-e^{i\beta_l}|}-4\log |\frac{y-Q_l^k}{y-e^{i\beta_l}}|)+E.
\end{equation}
Note that there is no summation of $l=0$ in (\ref{e-wk-important}).
Now for $s=1,...,N$, we consider the Pohozaev identity around $Q_s^k$. Let $\Omega_{s,k}=B(Q_s^k,r)$ for small $r>0$. For $v_k$ we have
\begin{align}\label{pi-vk}
&\int_{\Omega_{s,k}}\partial_{\xi}(|y|^{2N}\mathfrak{h}_k(\delta_ky))e^{v_k}-\int_{\partial\Omega_{s,k}}e^{v_k}|y|^{2N}\mathfrak{h}_k(\delta_ky)(\xi\cdot \nu)\\
=&\int_{\partial \Omega_{s,k}}(\partial_{\nu}v_k\partial_{\xi}v_k-\frac 12 |\nabla v_k|^2(\xi \cdot \nu))dS.\nonumber
\end{align}
where $\xi$ is an arbitrary unit vector. Correspondingly the Pohozaev identity for $V_k$ is
\begin{align}\label{pi-Vk}
&\int_{\Omega_{s,k}}\partial_{\xi}(|y|^{2N}\mathfrak{h}_k(\delta_ke_1))e^{V_k}-\int_{\partial\Omega_{s,k}}e^{V_k}|y|^{2N}\mathfrak{h}_k(\delta_ke_1)(\xi\cdot \nu)\\
=&\int_{\partial \Omega_{s,k}}(\partial_{\nu}V_k\partial_{\xi}V_k-\frac 12 |\nabla V_k|^2(\xi \cdot \nu))dS.\nonumber
\end{align}

 First we notice that the second term on the left hand side of (\ref{pi-vk}) or (\ref{pi-Vk}) is $O(e^{-\mu_k})$, so both terms are considered as errors. Now we first focus on (\ref{pi-vk}).
  After using the expansion for single equation as in the evaluation of (\ref{pi-e-1}), the first term on the left hand side of (\ref{pi-vk}) is:
\begin{align*}
&\int_{\Omega_{s,k}}\partial_{\xi}(|y|^{2N}\mathfrak{h}_k(\delta_ky))e^{v_k}\\
=&\int_{\Omega_{s,k}}\partial_{\xi}\bigg (\log (|y|^{2N}\mathfrak{h}_k(\delta_ky)) \bigg )|y|^{2N}\mathfrak{h}_k(\delta_ky)e^{v_k}\\
=&8\pi(2N\frac{Q_s^k}{|Q_s^k|^2}+\delta_k\nabla (\log \mathfrak{h}_k)(\delta_kQ_s^k))\cdot \xi+O(\mu_ke^{-\mu_k}).
\end{align*}
In a similar fashion, the first term of the left hand side of (\ref{pi-Vk}) is
$$\int_{\Omega_{s,k}}\partial_{\xi}(|y|^{2N}\mathfrak{h}_k(\delta_ke_1))e^{V_k}
=8\pi(2Ne^{i\beta_s}\cdot \xi)+O(\mu_ke^{-\mu_k}).
$$
For $l=1,..,N$, by $Q_l^k=e^{i\beta_l}(1+m_l^k)+O(\delta_k^2)$, we have
\begin{align*}
&\frac{Q_l^k}{|Q_l^k|^2}-e^{i\beta_l}
=\frac{1}{\bar Q_l^k}-\frac{1}{e^{-i\beta_l}}\\
=&\frac{e^{-i\beta_l}-\bar Q_l^k}{\bar Q_l^ke^{-i\beta_l}}
=-\frac{e^{-i\beta_l}\bar m_l^k}{e^{-2i\beta_l}}+E\\
=&-e^{i\beta_l}\bar m_l^k+E
=-\sum_ja^{lj}e^{i(\beta_l-\beta_j)}L_k\delta_k+E,
\end{align*}
where the last step is based on (\ref{loca-m}).
Thus the difference of the left hand sides gives the following leading term:
$$8\pi \delta_k\bigg ((1-2N\sum_{j=1}^Na^{sj}e^{i(\beta_s-\beta_j)})L_k \bigg )\cdot \xi+E. $$

To evaluate the right hand side, since $v_k=V_k+w_k+E$, the right hand side of (\ref{pi-vk}) is
\begin{align*}
&\int_{\partial\Omega_{s,k}}(\partial_{\nu} v_k\partial_{\xi} v_k-\frac 12 |\nabla v_k|^2(\xi\cdot \nu))dS\\
=&\int_{\partial \Omega_{s,k}}(\partial_{\nu} V_k\partial_{\xi} V_k-\frac 12 |\nabla  V_k|^2(\xi\cdot \nu))dS\\
&+\int_{\partial\Omega_{s,k}}(\partial_{\nu} V_k\partial_{\xi}w_k+\partial_{\nu}w_k\partial_{\xi} V_k-(\nabla  V_k\cdot \nabla w_k)(\xi \cdot \nu))dS+E.
\end{align*}
where we have used $w_k(y)=O(\delta_k)+E$.   Thus the difference of two pohozaev identities gives
\begin{align}\label{important-1}
&8\pi \delta_k(1-2N\sum_{j=1}^Na^{sj}e^{i(\beta_s-\beta_j)})\partial_{\xi}(\log \mathfrak{h}_k)(0)\\
=&\int_{\partial\Omega_{s,k}}(\partial_{\nu} V_k\partial_{\xi}w_k+\partial_{\nu}w_k\partial_{\xi} V_k-(\nabla  V_k\cdot \nabla w_k)(\xi \cdot \nu))dS+E\nonumber
\end{align}

Now we evaluate $\nabla w_k$ on $\partial B(Q_s^k,r)$ for $r>0$ fixed. For simplicity we omit $k$ in $m_i^k$ and $Q_l^k$.  For $y=e^{i\beta_s}+re^{i\theta}$ we obtain
from (\ref{e-wk-important}) that
\begin{align*}
\nabla w_k(y)&=\sum_{l=1}^N(-4\frac{y-Q_l}{|y-Q_l|^2}+4\frac{y-e^{i\beta_l}}{|y-e^{i\beta_l}|^2})+E\\
&=4\sum_{l=1}^N(\frac{1}{\bar y-e^{-i\beta_l}}-\frac{1}{\bar y-\bar Q_l})+E\\
&=4\sum_{l=1}^N\frac{e^{-i\beta_l}-\bar Q_l}{(\bar y-e^{-i\beta_l})(\bar y-\bar Q_l^k)}+E\\
&=4\sum_{l=1}^N\frac{-e^{-i\beta_l}\bar m_l}{(\bar y-e^{-i\beta_l})^2}+E.
\end{align*}
Separating $l=s$ with others, we have
$$
\nabla w_k(y)=-\frac{4}{r^2}e^{(2\theta-\beta_s)i}\bar m_s-4\sum_{l\neq s}^N\frac{e^{i\beta_l}\bar m_l}{(r e^{-i\theta}+e^{-i\beta_s}-e^{-i\beta_l})^2}+E.
$$

To evaluate the second term in $\nabla w_k$, we use these two formulas:
$$\sin \alpha-\sin \beta=2\cos \frac{\alpha+\beta}2\sin\frac{\alpha-\beta}2,\quad
\cos \alpha-\cos \beta=-2\sin\frac{\alpha+\beta}2\sin\frac{\alpha-\beta}2 $$
to obtain
\begin{align}\label{one-d}
&e^{-i\beta_s}-e^{-i\beta_l}=\cos\beta_s-\cos \beta_l+i(\sin \beta_l-\sin \beta_s)\nonumber\\
=&2\sin\frac{\beta_s+\beta_l}2\sin\frac{\beta_l-\beta_s}2+2i\cos \frac{\beta_s+\beta_l}2\sin\frac{\beta_l-\beta_s}2\nonumber\\
=&2\sin\frac{(l-s)\pi}{N+1}e^{i(\frac{\pi}2-\frac{l+s}{N+1}\pi)}.
\end{align}
Thus for $l\neq s$,
\begin{align*}
&\frac{e^{i\beta_l}\bar m_l}{(r e^{-i\theta}+e^{-i\beta_s}-e^{-i\beta_l})^2}
=\frac{e^{i\beta_l}\bar m_l}{(re^{-i\theta}+2\sin\frac{(l-s)\pi}{N+1}e^{i(\frac{\pi}2-\frac{l+s}{N+1}\pi)})^2}\\
=&\frac{e^{i\beta_l}\bar m_l}{4(\sin \frac{(l-s)\pi}{N+1})^2e^{i\pi-i\frac{2(l+s)}{N+1}\pi}(1+\frac{re^{-i\theta}}{2\sin\frac{(l-s)\pi}{N+1}e^{i(\frac{\pi}2-
\frac{l+s}{N+1}\pi)}})^2}
\end{align*}
Thus for $r$ small and fixed we have
\begin{align*}
&\frac{e^{i\beta_l}\bar m_l}{(r e^{-i\theta}+e^{-i\beta_s}-e^{-i\beta_l})^2}
=\frac{(-1)e^{i\beta_l}\bar m_le^{\frac{2(l+s)\pi}{N+1}i}}{4(\sin \frac{\pi(l-s)}{N+1})^2}
(1+O(r))\\
=&-\frac 14 d_{|l-s|}e^{\frac{4l+2s}{N+1}\pi i}\bar m_l+O(r).
\end{align*}
Thus we have the following important expression of $\nabla w_k$:
\begin{equation}\label{grad-w-k}
\nabla w_k(y)=-\frac{4}{r^2}e^{(2\theta-\beta_s)i}\bar m_s+\sum_{l\neq s}^Nd_{|l-s|} e^{\frac{4l+2s}{N+1}\pi i}\bar m_l+O(r)|L_k|\delta_k+E.
\end{equation}

On the other hand from the expression of $V_k$ in (\ref{Vk-exp}) we have, for $y=e^{i\beta_s}+re^{i\theta}$,
\begin{align*}
\nabla V_k(y)&=-4\sum_{l=1}^N\frac{y-e^{i\beta_l}}{|y-e^{i\beta_l}|^2}+O(\mu_ke^{-\mu_k})\\
&=(-4)\sum_{l=1}^N\frac{1}{re^{-i\theta}+e^{-i\beta_s}-e^{-i\beta_l}}+O(\mu_ke^{-\mu_k})\\
=&-\frac{4 e^{i\theta}}{r}-4\sum_{l\neq s}^N\frac{1}{re^{-i\theta}+e^{-i\beta_s}-e^{-i\beta_l}}+O(\mu_ke^{-\mu_k}).
\end{align*}
Here we note that $V_k$ has $N+1$ local maximums exactly located at $e^{i\frac{2\pi l}{N+1}}$ ($l=0,...,N$). The behavior of $V_k$ around each local maximum is very similar to that of a single Liouville equation with no singularity.
By similar computation in (\ref{one-d}) we have
\begin{align*}
re^{-i\theta}+e^{-i\beta_s}-e^{-i\beta_l}=r e^{-i\theta}+2\sin\frac{\pi(l-s)}{N+1} e^{i(\frac{\pi}2-\frac{l+s}{N+1}\pi)}\\
=2\sin\frac{\pi(l-s)}{N+1}e^{i(\frac{\pi}2-\frac{(l+s)\pi}{N+1})}(1+\frac{r e^{i(-\theta+\frac{(l+s)\pi}{N+1}-\frac{\pi}2)}}{2\sin \frac{\pi(l-s)}{N+1}})+O(r^2).
\end{align*}
Using this we can write $\nabla V_k$  as
\begin{align}\label{grad-v-k-r}
\nabla V_k=-\frac{4e^{i\theta}}r-4\sum_{l\neq s}^N\frac{e^{-i(\frac{\pi }2-\frac{(l+s)\pi}{N+1})}}{2\sin\frac{(l-s)\pi}{N+1}}
(1-\frac{ re^{i(-\theta+\frac{l+s}{N+1}\pi-\frac{\pi}2)}}{2\sin \frac{(l-s)\pi}{N+1}})+O(r^2)+E\nonumber \\
=-\frac{4e^{i\theta}}r+\sum_{l\neq s}^N\frac{2ie^{\frac{l+s}{N+1}\pi i}}{\sin \frac{\pi(l-s)}{N+1}}-\sum_{l\neq s}^N\frac{e^{\frac{2(l+s)\pi i}{N+1}-i\theta}}{\sin^2(\frac{\pi(l-s)}{N+1})}r+O(r^2)+E.
\end{align}
Correspondingly we have, by direct computation, that
\begin{equation}\label{d-nu-v}
\partial_{\nu}V_k=-\frac{4}{r}+\sum_{l\neq s}^N\frac{2\sin(\theta-\frac{(l+s)\pi}{N+1})}{\sin(\frac{\pi(l-s)}{N+1})}
-\sum_{l\neq s}^Nd_{|l-s|}\cos(\beta_{l+s}-2\theta)r+O(r^2)+E.
\end{equation}
For $\partial_{\xi}w_k$ we have, using $m_s=a_s+ib_s$ in (\ref{grad-w-k})
\begin{align*}
\partial_{\xi}w_k
=&-\frac{4}{r^2}\big ((a_s\cos(2\theta-\beta_s)+b_s\sin(2\theta-\beta_s))\xi_1\\
&+(a_s\sin (2\theta-\beta_s)-b_s\cos(2\theta-\beta_s))\xi_2\big )\\
&+\sum_{l\neq s}^Nd_{|l-s|}\big ((a_l\cos(\frac{4l+2s}{N+1}\pi)+b_l\sin(\frac{4l+2s}{N+1}\pi))\xi_1\\
&+(a_l\sin(\frac{4l+2s}{N+1}\pi)-b_l\cos(\frac{4l+2s}{N+1}\pi))\xi_2\big )
+O(r)|L_k|\delta_k+E.
\end{align*}

Then $\int_{\partial \Omega_{s,k}}\partial_{\nu}V_k\partial_{\xi}w_k$ comes from the integration of
\begin{align*}
&-\frac{4}r(\sum_{l\neq s}^Nd_{|l-s|}\big ( ( a_l \cos(\frac{4l+2s}{N+1}\pi)+b_l\sin(\frac{4l+2s}{N+1}\pi))\xi_1\\
&+(a_l \sin(\frac{4l+2s}{N+1}\pi)-b_l\cos(\frac{4l+2s}{N+1}\pi))\xi_2\big)\\
&+\frac 4r\sum_{l\neq s}^N d_{|l-s|}\cos(\beta_{l+s}-2\theta)\big ((a_s \cos(2\theta-\beta_s)+b_s\sin(2\theta-\beta_s))\xi_1\\
&+(a_s\sin(2\theta-\beta_s)-b_s \cos(2\theta-\beta_s))\xi_2\big ).
\end{align*}

For a part of the last term we have
\begin{align*}
&\cos(\beta_{l+s}-2\theta)\big ((a_s \cos(2\theta-\beta_s)+b_s\sin(2\theta-\beta_s))\xi_1\\
&+(a_s\sin(2\theta-\beta_s)-b_s \cos(2\theta-\beta_s))\xi_2\big )\\
=&(\cos \beta_{l+s}\cos 2\theta+\sin \beta_{l+s}\sin 2\theta)\\
&(a_s\cos 2\theta \cos \beta_s \xi_1+a_s\sin 2\theta \sin \beta_s  \xi_1+b_s \sin 2\theta \cos \beta_s \xi_1-b_s\cos 2\theta \sin \beta_s \xi_1)\\
&-b_s\cos 2\theta\cos \beta_s \xi_2 -b_s \sin 2\theta \sin \beta_s \xi_2 +a_s \sin 2\theta \cos \beta_s \xi_2-a_s \cos 2\theta \sin \beta_s \xi_2).
\end{align*}
Thus the integration of the last terms of $\partial_{\nu}V_k\partial_{\xi}w_k$ is
\begin{align*}
&2\pi r\int_{0}^{2\pi}\bigg (\frac 4r\sum_{l\neq s}^N d_{|l-s|}\cos(\beta_{l+s}-2\theta)\big ((a_s \cos(2\theta-\beta_s)+b_s\sin(2\theta-\beta_s))\xi_1\\
&+(a_s\sin(2\theta-\beta_s)-b_s \cos(2\theta-\beta_s))\xi_2\big )\bigg )d\theta\\
=&2\pi r\int_{0}^{2\pi}\bigg (\frac 4r\sum_{l\neq s}^N d_{|l-s|}\cos(\beta_{l+s}-2\theta)\big ((a_s \cos(2\theta-\beta_s)+b_s\sin(2\theta-\beta_s))\xi_1\\
&+(a_s\sin(2\theta-\beta_s)-b_s \cos(2\theta-\beta_s))\xi_2\big )\bigg )d\theta\\
=&4\pi\sum_{l\neq s}d_{|l-s|}\big ( \cos \beta_{l+s}(a_s \cos \beta_s \xi_1-b_s \sin \beta_s \xi_1-b_s \cos \beta_s \xi_2-a_s \sin \beta_s \xi_2)\\
&+ \sin \beta_{l+s}(a_s\sin \beta_s  \xi_1+b_s \cos \beta_s \xi_1-b_s \sin \beta_s \xi_2+a_s \cos \beta_s \xi_2) \big )\\
=&4\pi \sum_{l\neq s}d_{|l-s|}\big ( a_s (\cos \beta_l \xi_1+\sin \beta_l \xi_2)+\pi b_s (\sin \beta_l \xi_1-\cos \beta_l \xi_2)\big )\\
=&4\pi \sum_{l\neq s} d_{|l-s|}(\bar m_s e^{i\beta_l})\cdot \xi,
\end{align*}
and
\begin{equation}\label{one-part-pi}
\int_{\partial \Omega_{s,k}}\partial_{\nu}V_k\partial_{\xi}w_k=-8\pi\sum_{l\neq s}^Nd_{|l-s|}(e^{\frac{4l+2s}{N+1}\pi i}\bar m_l)\cdot \xi+4\pi\sum_{l\neq s}^Nd_{|l-s|}(\bar m_s
e^{i\beta_l})\cdot \xi+E.
\end{equation}
Next we compute the integral of $\partial_{\nu}w_k\partial_{\xi}V_k$. From the expression of $\nabla w_k$ we have
\begin{align*}
&\partial_{\nu}w_k
=-\frac{4}{r^2}\bigg ((a_s\cos(2\theta-\beta_s)+b_s\sin(2\theta-\beta_s))\cos \theta\\
&+(a_s\sin(2\theta-\beta_s)-b_s\cos (2\theta-\beta_s))\sin \theta\bigg )\\
&+\sum_{l\neq s}^Nd_{|l-s|}\bigg ( (a_l\cos(\frac{4l+2s}{N+1}\pi)+b_l\sin(\frac{4l+2s}{N+1}\pi))\cos \theta\\
&+(a_l\sin(\frac{4l+2s}{N+1}\pi)-b_l\cos(\frac{4l+2s}{N+1}\pi))\sin \theta \bigg )
+O(r)|L_k|\delta_k+E.
\end{align*}
After combining terms we have
\begin{align*}
\partial_{\nu}w_k&=-\frac{4}{r^2}(a_s\cos(\theta-\beta_s)+b_s\sin(\theta-\beta_s))\\
&+\sum_{l\neq s}^Nd_{|l-s|}(a_l\cos(\frac{4l+2s}{N+1}\pi-\theta)+b_l\sin(\frac{4l+2s}{N+1}\pi-\theta))+O(r)|L_k|\delta_k+E.
\end{align*}
From the expression of $\nabla V_k$ in (\ref{grad-v-k-r}) we have
\begin{align*}
\partial_{\xi}V_k&=-\frac{4}{r}(\cos \theta \xi_1+\sin \theta \xi_2)
+2\sum_{l\neq s}^N\frac{(-\sin(\frac{l+s}{N+1}\pi)\xi_1+\cos(\frac{l+s}{N+1}\pi)\xi_2)}{\sin\frac{\pi(l-s)}{N+1}}\\
&-\sum_{l\neq s}^Nd_{|l-s|}(\cos(\beta_{l+s}-\theta)\xi_1+\sin(\beta_{l+s}-\theta)\xi_2)r+E.
\end{align*}

First we consider the integral of
\begin{align*}
&(-\frac{4}{r^2})(a_s\cos(\theta-\beta_s)+b_s\sin(\theta-\beta_s))\\
&(-\frac 4r(\cos \theta \xi_1+\sin \theta \xi_2)-\sum_{l\neq s}^Nd_{|l-s|}(\cos(\beta_{l+s}-\theta)\xi_1+\sin(\beta_{l+s}-\theta)\xi_2)r)
\end{align*}
After applying $2\pi r\int_{0}^{2\pi} d\theta$ we have
\begin{align}\label{add-part-1}
&\frac{16\pi}{r^2}((a_s \cos \beta_s-b_s\sin \beta_s)\xi_1+(a_s\sin \beta_s +b_s\cos \beta_s)\xi_2)\nonumber\\
&+4\pi\sum_{l\neq s}^Nd_{|l-s|}(a_s\cos \beta_l \xi_1+a_s \sin \beta_l \xi_2+b_s\sin \beta_l \xi_1-b_s\cos \beta_l \xi_2)\nonumber\\
=&\frac{16\pi}{r^2}(m_se^{i\beta_s})\cdot \xi+4\pi\sum_{l\neq s}^Nd_{|l-s|}(\bar m_se^{i\beta_l})\cdot \xi.
\end{align}
This is one part of $\int_{\partial \Omega_{s,k}}\partial_{\nu}w_k\partial_{\xi}V_k$.  The other part comes from the integration of
$$\sum_{l\neq s}d_{|l-s|}(a_l\cos(\beta_{2l+s}-\theta)+b_l\sin(\beta_{2l+s}-\theta))(-\frac{4}{r})(\cos\theta \xi_1+\sin \theta \xi_2), $$
which is
\begin{align*}
&\sum_{l\neq s}^Nd_{|l-s|}\bigg ((a_l\cos(\frac{4l+2s}{N+1}\pi)+b_l\sin(\frac{4l+2s}{N+1}\pi))\cos \theta\\
&+(a_l\sin(\frac{4l+2s}{N+1}\pi)-b_l\cos(\frac{4l+2s}{N+1}\pi))\sin \theta\bigg )
(-\frac 4r)(\cos \theta \xi_1+\sin \theta \xi_2).
\end{align*}
After integration of the above we have
\begin{align}\label{add-part-2}
&-4\pi\sum_{l\neq s}^Nd_{|l-s|}\bigg (a_l\cos(\frac{4l+2s}{N+1}\pi)\xi_1+b_l\sin(\frac{4l+2s}{N+1}\pi)\xi_1 \nonumber \\
&+a_l \sin(\frac{4l+2s}{N+1}\pi)\xi_2-b_l\cos(\frac{4l+2s}{N+1}\pi)\xi_2\bigg )\nonumber\\
=&-4\pi \sum_{l\neq s}^Nd_{|l-s|}(\bar m_le^{i\frac{4l+2s}{N+1}\pi})\cdot \xi.
\end{align}
Putting (\ref{add-part-1}) and (\ref{add-part-2}) together we have
\begin{equation}\label{add-3}
\int_{\partial \Omega_{s,k}}\partial_{\nu}w_k\partial_{\xi}V_k=
\bigg (\frac{16\pi}{r^2}m_se^{i\beta_s}+4\pi\sum_{l\neq s}^N d_{|l-s|}(\bar m_se^{i\beta_l}-\bar m_le^{i\frac{4l+2s}{N+1}\pi})\bigg )\cdot \xi+E.
\end{equation}

Now we compute $\int_{\partial \Omega_{s,k}}(\nabla V_k\cdot \nabla w_k)(\xi\cdot \nu)$. By the expressions of $\nabla V_k$ and $\nabla w_k$, the terms of $(\nabla V_k\cdot \nabla w_k)(\xi\cdot \nu)$ that survive integration are
\begin{align*}
&[\frac{16}{r^3}(e^{(2\theta-\beta_s)i}\bar m_s)\cdot e^{i\theta}+\frac 4r(e^{i(2\theta-\beta_s)}\bar m_s)\cdot (\sum_{l\neq s}^Nd_{|l-s|}e^{i(\beta_{l+s}-\theta)})\\
&+\sum_{l\neq s}^Nd_{|l-s|}(e^{\frac{4l+2s}{N+1}\pi i}\bar m_l)\cdot (-\frac{4e^{i\theta}}r)](\xi_1\cos \theta+\xi_2\sin \theta).\\
=&[\frac{16}{r^3}((a_s\cos (2\theta-\beta_s)+b_s\sin (2\theta-\beta_s))\cos \theta\\
&+(a_s\sin (2\theta-\beta_s)-b_s\cos (2\theta-\beta_s))\sin \theta)\\
&+\frac{4}{r}\sum_{l\neq s}^Nd_{|l-s|}((a_s\cos(2\theta-\beta_s)+b_s\sin(2\theta-\beta_s))\cos(\beta_{l+s}-\theta)\\
&+(a_s\sin(2\theta-\beta_s)-b_s\cos(2\theta-\beta_s))\sin(\beta_{l+s}-\theta))\\
&-\frac{4}{r}\sum_{l\neq s}^Nd_{|l-s|}((a_l\cos(\frac{4l+2s}{N+1}\pi)+b_l\sin(\frac{4l+2s}{N+1}\pi))\cos\theta\\
&+(a_l\sin(\frac{4l+2s}{N+1}\pi)-b_l\cos (\frac{4l+2s}{N+1}\pi))\sin \theta)]
(\xi_1\cos \theta+\xi_2 \sin \theta).
\end{align*}

After simplification these terms become
\begin{align*}
&[\frac{16}{r^3}(a_s\cos(\theta-\beta_s)+b_s\sin(\theta-\beta_s))\\
&+\frac{4}{r}\sum_{l\neq s} d_{|l-s|}(a_s\cos(3\theta-\frac{2l+4s}{N+1}\pi)+b_s\sin (3\theta-\frac{2l+4s}{N+1}\pi))\\
&-\frac{4}{r}\sum_{l\neq s} d_{|l-s|}(a_l\cos(\frac{4l+2s}{N+1}\pi-\theta)+b_l\sin (\frac{4l+2s}{N+1}\pi-\theta))](\xi_1 \cos \theta+\xi_2 \sin \theta).
\end{align*}
After integration we have
\begin{align}\label{add-4}
&\int_{\partial \Omega_{s,k}}(\nabla V_k\cdot \nabla w_k)(\xi\cdot \nu)\\
&=\frac{16\pi}{r^2}(a_s\cos \beta_s\xi_1+a_s\sin \beta_s\xi_2+b_s\xi_2\cos \beta_s-b_s\sin \beta_s\xi_1)\nonumber\\
&-4\pi\sum_{l\neq s}d_{|l-s|}(a_l\cos(\frac{4l+2s}{N+1}\pi)\xi_1+a_l\sin (\frac{4l+2s}{N+1}\pi)\xi_2\nonumber\\
&+b_l \sin (\frac{4l+2s}{N+1}\pi)\xi_1
-b_l\cos (\frac{4l+2s}{N+1}\pi)\xi_2)+E\nonumber\\
&=\frac{16\pi}{r^2}(m_se^{i\beta_s})\cdot \xi
-4\pi\sum_{l\neq s}d_{|l-s|}(\bar m_l e^{i\frac{4l+2s}{N+1}\pi})\cdot \xi+E.\nonumber
\end{align}
Combining (\ref{one-part-pi}), (\ref{add-3}) and (\ref{add-4}) we have
\begin{align}\label{2-part-pi}
\int_{\partial \Omega_{s,k}}(\partial_{\nu}V_k\partial_{\xi}w_k+\partial_{\xi}V_k\partial_{\nu}w_k-(\nabla V_k\cdot \nabla w_k)(\xi\cdot \nu)  \\
=8\pi\sum_{l\neq s}d_{|l-s|}(\bar m_se^{i\beta_l}-\bar m_le^{i\beta_{2l+s}})\cdot \xi+E. \nonumber
\end{align}
Using (\ref{loca-m}) and (\ref{2-part-pi}) in (\ref{important-1})
we see that as long as one of the following equations is violated
\begin{align}\label{cru-coe}
&1-2N\sum_{j=1}^Na^{sj}e^{i(\beta_s-\beta_j)}\\
=&\sum_{l=1,l\neq s}^Nd_{|l-s|}\sum_{j=1}^N(a^{sj}e^{i(\beta_l-\beta_j)}-a^{lj}e^{i(\beta_{2l+s}-\beta_j)}), \quad
s=1,..,N \nonumber
\end{align}
 the vanishing rate is proved. But here observe that the vanishing rate of $\nabla \log \mathfrak{h}_k(0)$ holds for $N=1$, as the coefficient on the left is $-1$ and is $0$ on the right. In order to prove the general case for $N\ge 2$ we assume that (\ref{cru-coe}) holds for all $s=1,...,N$ and will derive a contradiction later. First we establish a few lemmas to be used later.

\begin{lem}\label{compu-1}
$$D-d_1e^{i\beta_1}-d_2e^{i\beta_2}-...-d_Ne^{i\beta_N}=2N. $$
\end{lem}
\noindent{\bf Proof:} If $N$ is odd, we observe that $d_{(N+1)/2}=1$, $e^{i\beta_{(N+1)/2}}=-1$, using $d_l=d_{N+1-l}$ and $e^{i\beta_l}=e^{-i\beta_{N+1-l}}$ we have
\begin{align*}
&D-(d_1e^{i\beta_1}+...+d_Ne^{i\beta_N})\\
=&D-2\sum_{j=1}^{\frac{N-1}2}\frac{\cos \frac{2\pi j}{N+1}}{\sin^2(\frac{j\pi}{N+1})}+1=D-\sum_{j=1}^{\frac{N-1}2}\frac{2\cos(\frac{2\pi j}{N+1})-2+2}{\sin^2(\frac{j\pi}{N+1})}+1\\
=&D+\sum_{j=1}^{\frac{N-1}2}\frac{4\sin^2(\frac{j\pi}{N+1})}{\sin^2(\frac{j\pi}{N+1})}-2\sum_{j=1}^{\frac{N-1}{2}}\frac{1}{\sin^2(\frac{j\pi}{N+1})}+1\\
=&D+4(\frac{N-1}2)-(D-1)+1=2N.
\end{align*}
If $N$ is even,
\begin{align*}
&d_1e^{i\beta_1}+...+d_Ne^{i\beta_N}=\sum_{l=1}^{N/2}d_le^{i\beta_l}+\sum_{l=N/2+1}^Nd_le^{i\beta_l}\\
=&2\sum_{l=1}^{N/2}\frac{\cos(\frac{2l\pi}{N+1})-1+1}{\sin^2(\frac{l\pi}{N+1})}=2\sum_{l=1}^{N/2}\frac{1-2\sin^2(\frac{l\pi}{N+1})}{\sin^2(\frac{l\pi}{N+1})}
=-2N+D.
\end{align*}
Lemma \ref{compu-1} is established. $\Box $

\begin{lem}\label{compu-8} Let $A=(a^{ij})_{N\times N}$, then
\begin{equation}\label{imp-3}
\sum_{j=1}^Na^{ij}d_{N+1-j}=1,~i=1,\cdots,N,
\end{equation}
\end{lem}

\noindent{\bf Proof of Lemma \ref{compu-8}:}
Reading the first row of $A^{-1}A=I$ we have the following:
\begin{equation*}
\begin{aligned}
a^{11}D-a^{12}d_1-\cdots -a^{1N}d_{N-1}=1,\\
-a^{11}d_1+a^{12}D-\cdots -a^{1N}d_{N-2}=0,\\
\vdots\\
-a^{11}d_{N-1}-a^{12}d_{N-2}-\cdots +a^{1N}D=0,\\
\end{aligned}
\end{equation*}
adding them together, it gives
\begin{equation*}
a^{11}d_N+a^{12}d_{N-1}+\cdots +a^{1N}d_1=1,
\end{equation*}
similarly, (\ref{imp-3}) can be obtained by the same analysis on other rows. Lemma \ref{compu-8} is established. $\Box$

\begin{lem}\label{compu-2}
$$\sum_{s,t}a^{st}e^{i\beta_s}e^{-i\beta_t}=\frac{N+1}{2N}.$$
\end{lem}

\noindent{\bf Proof:} We set $B$ to be the circular matrix
$$
B=\left[\begin{matrix}
D&-d_1&\cdots&-d_{N}\\
-d_N&D&\cdots&-d_{N-1}\\
\vdots&\vdots&\ddots&\vdots\\
-d_1&-d_2&\cdots&D
\end{matrix}\right]=
\left(\begin{matrix}
D& -\vec{d}\\
-\vec{d}^t & A
\end{matrix}
\right),
$$
where $\vec{d}=(d_1,\cdots,d_N).$  Note that we have used $d_l=d_{N+1-l}$. Using the definition of $\beta_l$ we have
$$
B\left(\begin{matrix}
1\\e^{i\beta_1}\\ \vdots \\ e^{i\beta_N}
\end{matrix}\right)
=(D-d_1e^{i\beta_1}-\cdots-d_Ne^{i\beta_N})\left(\begin{matrix}
1\\e^{i\beta_1}\\ \vdots \\ e^{i\beta_N}
\end{matrix}\right),
$$
Comparing the last $N$ rows, we have, by Lemma \ref{compu-1} that
\begin{equation}\label{imp-1}
-\left(\begin{matrix}
d_N\\ d_{N-1}\\ \vdots \\ d_1
\end{matrix}\right)+A \left(\begin{matrix}
e^{i\beta_1}\\e^{i\beta_2}\\ \vdots \\ e^{i\beta_N}
\end{matrix}\right)
=2N\left(\begin{matrix}
e^{i\beta_1}\\e^{i\beta_2}\\ \vdots \\ e^{i\beta_N}
\end{matrix}\right),
\end{equation}
Multiplying $(e^{i\beta_1},e^{i\beta_2},\cdots,e^{i\beta_N})A^{-1}$ to the conjugate of the above, we get
\begin{equation}
\label{compu-3}
\begin{aligned}
&-(e^{i\beta_1},e^{i\beta_2},\cdots,e^{i\beta_N})A^{-1}\left(\begin{matrix}
d_N\\ d_{N-1}\\ \vdots \\ d_1
\end{matrix}\right)+(e^{i\beta_1},e^{i\beta_2},\cdots,e^{i\beta_N})\left(\begin{matrix}
e^{-\beta_1}\\ e^{-i\beta_2}\\ \vdots \\ e^{-i\beta_N}
\end{matrix}\right)\\
&=2N(e^{i\beta_1}, e^{i\beta_2},\cdots,e^{i\beta_N})A^{-1}\left(\begin{matrix}
e^{-i\beta_1}\\ e^{-i\beta_2}\\ \vdots \\ e^{-i\beta_N}
\end{matrix}\right)
\end{aligned}
\end{equation}
The left hand side of \eqref{compu-3}:
we obtain from (\ref{imp-3}) that
\begin{align}\label{imp-2}
&-(e^{i\beta_1},e^{i\beta_2},\cdots,e^{i\beta_N})A^{-1}\left(\begin{matrix}
d_N\\ d_{N-1}\\ \vdots \\ d_1
\end{matrix}\right)\\
=&-(e^{i\beta_1},e^{i\beta_2},\cdots,e^{i\beta_N})\left(\begin{matrix}
1\\1\\ \vdots \\ 1
\end{matrix}\right)=-e^{i\beta_1}-e^{i\beta_2}-\cdots-e^{i\beta_N}=1. \nonumber
\end{align}
Lemma \ref{compu-2} is established after using (\ref{imp-2}) in (\ref{compu-3}). $\Box$

\medskip
\begin{lem}\label{add-lem}
$$\sum_{s,t=1}^ne^{i\beta_s}a^{st}e^{i\beta_t}=0.$$
\end{lem}
\noindent{\bf Proof of Lemma \ref{add-lem}:}
We multiply $(e^{i\beta_1},e^{i\beta_2},\cdots,e^{i\beta_N})A^{-1}$ from the left to (\ref{imp-1}) and obtain
\begin{equation}
\label{compu-10}
\begin{aligned}
&-(e^{i\beta_1},e^{i\beta_2},\cdots,e^{i\beta_N})A^{-1}\left(\begin{matrix}
d_N\\ d_{N-1}\\ \vdots \\ d_1
\end{matrix}\right)+\sum_{l=1}^Ne^{2i\beta_l}\left(\begin{matrix}
e^{-\beta_1}\\ e^{-i\beta_2}\\ \vdots \\ e^{-i\beta_N}
\end{matrix}\right)\\
&=2N(e^{i\beta_1}, e^{i\beta_2},\cdots,e^{i\beta_N})A^{-1}\left(\begin{matrix}
e^{i\beta_1}\\ e^{i\beta_2}\\ \vdots \\ e^{i\beta_N}
\end{matrix}\right)
\end{aligned}
\end{equation}
The first term on the left is $1$, while the second term is $-1$. Thus
Lemma \ref{add-lem} is established. $\Box$

\medskip

\begin{lem}\label{inv-lem}
 For each $s=1,...,N$,
 $$\sum_{j=1}^Na^{sj}e^{i\beta_j}=\frac{e^{i\beta_s}-1}{2N},\quad \sum_{j=1}^Na^{sj}e^{-i\beta_j}=\frac{e^{-i\beta_s}-1}{2N}. $$
\end{lem}

\noindent{\bf Proof of Lemma \ref{inv-lem}:}
Multiplying $A^{-1}$ to both sides of (\ref{imp-1}) we have
$$-A^{-1}\left(\begin{array}{c}
d_{N}\\
d_{N-1}\\
\vdots\\
d_1
\end{array}
\right)+
\left(\begin{array}{c}
e^{i\beta_1}\\
e^{i\beta_2}\\
\vdots\\
e^{i\beta_N}
\end{array}
\right)=2N A^{-1}\left(\begin{array}{c}
e^{i\beta_1}\\
e^{i\beta_2}\\
\vdots\\
e^{i\beta_N}
\end{array}
\right)$$
The first term on the left is evaluated by (\ref{imp-3}), thus Lemma \ref{inv-lem} follows immediately. $\Box$

\medskip

Based on Lemma \ref{inv-lem} we take the summation of $s$ for both sides of (\ref{cru-coe}). The left hand side is
\begin{align}\label{add-left-1}
&\sum_{s=1}^N(1-2N\sum_{j=1}^Na^{sj}e^{-i\beta_j}e^{i\beta_s})\\
=&N-\sum_{s=1}^N(e^{-i\beta_s}-1)e^{i\beta_s}=N-N+(-1)=-1.\nonumber
\end{align}
Now we evaluate the two terms on the right hands side of (\ref{cru-coe}). The first term is
\begin{equation}\label{fin-term-1}
\sum_{s=1}^N\sum_{l=1,l\neq s}^Nd_{|l-s|}e^{i\beta_l}\sum_{j=1}^Na^{sj}e^{-i\beta_j}
=\frac{1}{2N}\sum_{s=1}^N\sum_{l=1,l\neq s}^Nd_{|l-s|}e^{i\beta_l}(e^{-i\beta_s}-1).
\end{equation}
and the second term is
\begin{equation}\label{fin-term-2}
\sum_{s=1}^N\sum_{l\neq s}d_{|l-s|}\sum_{j=1}^Na^{lj}e^{i(\beta_{2l+s}-\beta_j)}=\frac{1}{2N}
\sum_{s=1}^N\sum_{l\neq s}d_{|l-s|}e^{i\beta_{2l+s}}(e^{-i\beta_l}-1).
\end{equation}
where Lemma \ref{inv-lem} is used.
First we evaluate (\ref{fin-term-1}),
let
$$\Lambda=\sum_{l=1}^Nd_le^{i\beta_l}=D-2N, $$
then
\begin{align}\label{add-11}
&\sum_{s=1}^N\sum_{l=1,l\neq s}^Nd_{|l-s|}e^{i\beta_l-i\beta_s}=(d_1e^{i\beta_1}+...+d_{N-1}e^{i\beta_{N-1}})\nonumber\\
&+(d_1e^{-i\beta_1}+d_1e^{i\beta_1}+d_2e^{i\beta_2}+...+d_{N-2}e^{i\beta_{N-2}})\nonumber\\
&+...+(d_{N-1}e^{-i\beta_{N-1}}+d_{N-2}e^{-i\beta_{N-2}}+...+d_{1}e^{-i\beta_1})\nonumber\\
&=(\Lambda-d_1e^{-i\beta_1})+(\Lambda-d_2e^{-i\beta_2})+...+(\Lambda-d_Ne^{-i\beta_N})\nonumber\\
&=(N-1)\Lambda.
\end{align}

\begin{align}\label{add-12}
&\sum_{s=1}^N\sum_{l=1,l\neq s}^Nd_{|l-s|}e^{i\beta_l}=\sum_{s=1}^Ne^{i\beta_s}\sum_{l\neq s}d_{|l-s|}e^{i(\beta_l-\beta_s)}\nonumber\\
=&e^{i\beta_1}(\Lambda-d_1e^{-i\beta_1})+e^{i\beta_2}(\Lambda-d_2e^{-i\beta_2})+...+e^{i\beta_N}(\Lambda-d_Ne^{-i\beta_N})\nonumber\\
=&-\Lambda-D.
\end{align}
Putting (\ref{add-11}) and (\ref{add-12}) together we have
\begin{equation}\label{tba-1}
\sum_{s=1}^N\sum_{l=1,l\neq s}^Nd_{|l-s|}e^{i\beta_l}\sum_{j=1}^Na^{sj}e^{-i\beta_j}=\frac{N\Lambda+D}{2N}.
\end{equation}
Then we evaluate (\ref{fin-term-2}):
\begin{equation}\label{add-13}
\frac{1}{2N}\bigg (\sum_{s=1}^N\sum_{l=1,l\neq s}^Nd_{|l-s|}e^{i\beta_{2l}}e^{i\beta_{s-l}}-\sum_{s=1}^N\sum_{l=1,l\neq s}^Nd_{|l-s|}e^{i\beta_{3l}}e^{i\beta_{s-l}}\bigg ).
\end{equation}
The first term in the brackets of (\ref{add-13}) is equal to
\begin{align}\label{add-15}
&\sum_s\sum_{l\neq s}d_{|l-s|}e^{i\beta_{2l}}e^{i\beta_{s-l}}=e^{i\beta_2}(d_1e^{i\beta_1}+d_2e^{i\beta_2}+...+d_{N-1}e^{i\beta_{N-1}})\nonumber\\
&+e^{i\beta_4}(d_1e^{-i\beta_1}+d_1e^{i\beta_1}+...+d_{N-2}e^{i\beta_{N-2}})\nonumber\\
&+...+e^{i\beta_{2N}}(d_{N-1}e^{-i\beta_{N-1}}+d_{N-2}e^{-i\beta_{N-2}}+...+d_1e^{-i\beta_1})\nonumber\\
=&e^{i\beta_2}(\Lambda-d_1e^{-i\beta_1})+e^{i\beta_4}(\Lambda-d_2e^{-i\beta_2})+...+e^{i\beta_{2N}}(\Lambda-d_Ne^{-i\beta_N})\nonumber\\
=&\Lambda(e^{i\beta_2}+e^{i\beta_4}+...+e^{i\beta_{2N}})-d_1e^{i\beta_1}-...-d_Ne^{i\beta_N}\nonumber\\
=&-2\Lambda.
\end{align}
For the second term in the brackets of (\ref{add-13}) we have
\begin{align}\label{add-14}
&\sum_s\sum_{l\neq s}d_{|l-s|}e^{i\beta_{3l}}e^{i\beta_{s-l}}=e^{i\beta_3}(d_1e^{i\beta_1}+d_2e^{i\beta_2}+...+d_{N-1}e^{i\beta_{N-1}})\nonumber\\
&+e^{i\beta_6}(d_1e^{-i\beta_1}+d_1e^{i\beta_1}+...+d_{N-2}e^{i\beta_{N-2}})\nonumber\\
&+...+e^{i\beta_{3N}}(d_{N-1}e^{-\beta_{N-1}}+d_{N-2}e^{-i\beta_{N-2}}+...+d_1e^{-\beta_1})\nonumber\\
=&e^{i\beta_3}(\Lambda-d_1e^{-i\beta_1})+e^{i\beta_6}(\Lambda-d_2e^{-i\beta_2})+...+e^{i\beta_{3N}}(\Lambda-d_Ne^{-i\beta_N})\nonumber\\
=&\Lambda(e^{i\beta_3}+e^{i\beta_6}+...+e^{i\beta_{3N}})-d_1e^{i\beta_2}-...-d_Ne^{i\beta_{2N}}\nonumber\\
=&-\Lambda-d_1e^{i\beta_2}-...-d_Ne^{i\beta_{2N}}.
\end{align}
The evaluation of $\sum_{l=1}d_le^{i\beta_{2l}}$ is divided into two cases depending on $N$ even or odd.
If $N$ is odd
\begin{align*}
&\sum_{l=1}^Nd_le^{i\frac{4l\pi}{N+1}}=\sum_{l=1}^{(N-1)/2}d_le^{i\frac{4l\pi}{N+1}}+1+\sum_{l=\frac{N+3}2}^Nd_le^{i\frac{4l\pi}{N+1}}\\
=&2\sum_{l=1}^{(N-1)/2}d_l\cos(\frac{4l\pi}{N+1})+1=(-4)\sum_{l=1}^{\frac{N-1}2}\frac{\sin^2(\frac{2l\pi}{N+1})}{\sin^2(\frac{l\pi}{N+1})}+D-1+1\\
=&-16\sum_{l=1}^{(N-1)/2}\cos^2(\frac{2l\pi}{N+1})+D=(-8)\sum_{l=1}^{(N-1)/2}(1+\cos(\frac{2l\pi}{N+1}))+D\\
=&-4(N-1)+D.
\end{align*}
If $N$ is even
\begin{align*}
&\sum_{l=1}^Nd_le^{i\frac{4l\pi}{N+1}}=\sum_{l=1}^{N/2}d_le^{i\frac{4l\pi}{N+1}}+\sum_{l=\frac N2+1}^Nd_le^{i\frac{4l\pi}{N+1}}\\
=&2\sum_{l=1}^{N/2}d_l\cos(\frac{4l\pi}{N+1})=2\sum_{l=1}^{N/2}d_l(\cos(\frac{4l\pi}{N+1})-1)+D\\
=&(-4)\sum_{l=1}^{\frac N2}\frac{\sin^2(\frac{2l\pi}{N+1})}{\sin^2(\frac{l\pi}{N+1})}+D
=-16\sum_{l=1}^{N/2}\cos^2(\frac{2l\pi}{N+1})+D\\
=&(-8)\sum_{l=1}^{N/2}(1+\cos(\frac{2l\pi}{N+1}))+D=-4N-8\sum_{l=1}^{N/2}\cos(\frac{2l\pi}{N+1})+D\\
=&-4N-8\sum_{l=1}^{N/2}\frac{\sin \frac{\pi}{N+1} \cos (\frac{2l\pi}{N+1})}{\sin \frac{\pi}{N+1}}+D\\
=&-4N-8\sum_{l=1}^{\frac{N}2}\frac{\sin\frac{(2l+1)\pi}{N+1}-\sin\frac{(2l-1)\pi}{N+1}}{2\sin\frac{\pi}{N+1}}+D=-4(N-1)+D.
\end{align*}

Thus in either case $\sum_{l=1}^Nd_le^{i\beta_{2l}}=-4(N-1)+D$, then the result of (\ref{add-14}) is $-\Lambda+4(N-1)-D$. Using this and (\ref{add-15}) in (\ref{add-13}) we have
\begin{align}\label{imp-9}
&\sum_{s=1}^N\sum_{l\neq s}d_{|l-s|}\sum_{j=1}^Na^{lj}e^{i(\beta_{2l+s}-\beta_j)}\\
=&\frac{1}{2N}(-2\Lambda+\Lambda-4(N-1)+D)=\frac 2N-1. \nonumber
\end{align}

Thus the combination of the two terms on the right hand side of (\ref{cru-coe}) is
$$\frac{\Lambda}2+\frac{D}{2N}+1-\frac 2N\neq -1,\quad N>1, N\in \mathbb N. $$
Thus Theorem \ref{Van-n2} is established. $\Box$

\section{Approximate $u_k$ by global solutions}

First we note that for simple blowup solutions, the approximation of $u_k$ using global solutions is much easier.  This part will be discussed in the appendix.

\bigskip

\noindent{\bf Proof of Theorem \ref{main-thm-2}:}   The assumption is $\delta_ke^{\mu_k/4}\le c_0$. Fixing the neighborhood of one $Q_m^k$, the expansion of $v_k$ is, taking $Q_m^k$ as the origin,
\begin{equation}\label{tem-vk-2}
v_k(y)=\log \frac{e^{\mu_{k,m}}}{(1+e^{\mu_{k,m}}\frac{|\tilde Q_m^k|^{2N}\mathfrak{h}_k(\delta_k\tilde Q_m^k)}8|y-\tilde Q_m^k|^2)^2}+\phi_m^k(y)+O(\mu_k^2e^{-\mu_k})
\end{equation}
where $\mu_{k,m}= v_k(\tilde Q_m^k)$. First we claim that
\begin{equation}\label{mu-lk}
\mu_{k,m}-\mu_k=O(\delta_k)+O(\mu_k^2e^{-\mu_k}).
\end{equation}
From the  Green's representation formula for $v_k$, we have, for $y$ away from bubbling areas and $|y|\sim 1$,
\begin{align*}
v_k(y)&=v_k|_{\partial \Omega_k}+\int_{\Omega_k}G(y,\eta)\mathfrak{h}_k(\eta)|\eta |^{2N}e^{v_k}d\eta, \\
&=v_k|_{\partial \Omega_k}+\sum_{l=0}^NG(y,Q_l^k)\int_{B(Q_l^k,\epsilon)}|\eta |^{2N}\mathfrak{h}_k(\delta_k\eta)e^{v_k}d\eta \\ &+\sum_{l}\int_{\Omega_k}(G(y,\eta)-G(y,Q_l^k))|\eta|^{2N}\mathfrak{h}_k(\delta_k\eta)e^{v_k}d\eta+O(e^{-\mu_k}),\\
&=v_k|_{\partial \Omega_k}+8\pi \sum_l G(y,Q_l^k)+E
\end{align*}
where $\Omega_k=B(0,\tau \delta_k^{-1})$.
 In particular if we consider $y$ located at $|y-Q_m^k|=\frac{\epsilon}2$, the expression of $v_k$ can be written as
\begin{align}\label{tem-vk}
v_k(y)&=v_k|_{\partial \Omega_k}-4\log |y-Q_m^k|+\phi_m^k\\
&-4\sum_{l=0,l\neq m}^N\log |Q_m^k-Q_l^k|+8\pi\sum_{l=0}^N H(Q_m^k,Q_l^k)+E,\nonumber
\end{align}
where
$$\phi_m^k=\sum_{l=0,l\neq m}^N(-4)\log \frac{|y-Q_l^k|}{|Q_m^k-Q_l^k|}+8\pi\sum_{l=0}^N(H(y,Q_l^k)-H(Q_m^k,Q_l^k)). $$
Comparing (\ref{tem-vk}) and (\ref{tem-vk-2}) we have
\begin{align}\label{tem-vk-3}
&-\mu_{m,k}-\log \frac{|\tilde Q_m^k|^{2N}\mathfrak{h}_k(\delta_k\tilde Q_m^k)}{8}\\
=&-4\sum_{l=0,l\neq m}^N\log |Q_m^k-Q_l^k|+8\pi\sum_{l=0}^N H(Q_m^k,Q_l^k)+v_k|_{\partial \Omega_k}+O(\mu_k^2e^{-\mu_k}). \nonumber
\end{align}
To evaluate terms in (\ref{tem-vk-3}) we observe that
$$|\tilde Q_m^k|^{2N}=1+E,\qquad \mathfrak{h}_k(\delta_k\tilde Q_m^k)=1+E,$$
$$Q_m^k=e^{\frac{2\pi m}{N+1}i}+E, \qquad \tilde Q_m^k=Q_m^k+O(e^{-\mu_k}),$$
and by  the expression of $H(y,\eta)$ we have
$$H(Q_m^k,Q_l^k)=\frac{1}{2\pi}\log (\tau \delta_k^{-1})+E.$$
Thus two terms in (\ref{tem-vk-3}) are
\begin{equation}\label{tri-1}
8\pi\sum_{l=0}^N H(Q_m^k,Q_l^k)=4(N+1)\log (\tau \delta_k^{-1})+E
\end{equation}
\begin{align}\label{tri-2}
\sum_{l=0,l\neq m}^N\log |Q_m^k-Q_l^k|&=\sum_{l=0,l\neq m}^N\log |e^{\frac{2\pi mi}{N+1}}-e^{\frac{2\pi li}{N+1}}|+E\\
&=\log (N+1)+E. \nonumber
\end{align}

Using (\ref{tri-1}) and (\ref{tri-2}) in (\ref{tem-vk-3}) we have
\begin{align}\label{vk-bry}
v_k|_{\partial \Omega_k}=&-\mu_{m,k}+\log 8+4\log (1+N)-4(1+N)\log (\tau \delta_k^{-1})\\
&+O(\delta_k^2)+O(\mu_k^2e^{-\mu_k}),\quad m=0,1,...,N. \nonumber
\end{align}
Clearly from (\ref{vk-bry}) we see that (\ref{mu-lk}) holds.
In order to approximate $v_k$ with a global solution we find $U_k$ which exactly has local maximums located at $e^{\frac{2\pi l}{N+1}i}$ and $U_k(e_1)=\mu_k$:
$$U_k(x)=\log \frac{e^{\mu_k}}{(1+\frac{\mathfrak{h}_k(\delta_ke_1)e^{\mu_k}}{8(1+N)^2}|y^{N+1}-e_1|^2)^2}, $$
where $e_1=(1,0)$ on $\mathbb R^2$.

First in the region $B(Q_l^k, e^{-\mu_k/2})$, the comparison between $v_k$ and $U_k$ boils down to the evaluation of:

\begin{equation}\label{glob-a}
\log \frac{e^{\mu_{l,k}}}{(1+\frac{\mathfrak{h}_k(\delta_ke_1)e^{\mu_{l,k}}}8|y-p_k|^2)^2}-\log \frac{e^{\mu_k}}{(1+\frac{\mathfrak{h}_k(\delta_ke_1)e^{\mu_k}}8|y|^2)^2},
\end{equation}
for $|p_k|=E$. By elementary computation we see that the difference between the two terms in (\ref{glob-a}) is $O(\delta_k^2e^{\mu_k/2})+O(\mu_ke^{-\mu_k/2})$ if $|y|\le C e^{-\mu_k/2}$. On the other hand, for $Ce^{-\mu_k/2}<|y|<\epsilon/2$, the comparison of expressions of $v_k$ and $U_k$ leads to the same conclusion. Moreover
$$v_k-U_k=O(\delta_k^2)+O(\mu_k^2e^{-\mu_k})\quad \mbox{on}\quad \partial B(Q_l^k,\epsilon).$$
Also we observe from the expression of $U_k$ that
$$v_k-U_k=O(\delta_k^2)+O(\mu_k^2e^{-\mu_k})\quad \mbox{on}\quad \partial \Omega_k. $$
Thus we obtain
the closeness of $v_k$ and $U_k$ on $\Omega_k\setminus (\cup_lB(Q_l^k,\epsilon/2))$  by the smallness of $v_k-U_k$ on $\partial \Omega_k$ and standard estimates by Green's representation formula.  Theorem \ref{main-thm-2} is established. $\Box$

\section{Vanishing of Second order derivatives}

The main purpose of this section is to prove the vanishing rate of the second derivatives of $\log \mathfrak{h}_k(0)$ in Theorem \ref{Van-n1}.

\medskip

\noindent{\bf Proof of Theorem \ref{Van-n1}:}  Obviously if $\delta_k^2\le C\mu_ke^{-\mu_k}$ for some $C>0$ there is nothing to prove. So in this section we assume
$\mu_ke^{-\mu_k}=o(1)\delta_k^2$.

Observing (\ref{pi-each-p}) and (\ref{gra-each-3}) we see that the equation of Pohozaev identity now becomes
\begin{equation}\label{large-delta-1}
2N\frac{Q_l^k}{|Q_l^k|^2}-4\sum_{m\neq l}\frac{Q^k_l-Q^k_m}{|Q^k_l-Q^k_m|^2}+\delta_k \nabla (\log \mathfrak{h}_k)(\delta_k Q^k_l)=O(\delta_k^3)+O(\mu_ke^{-\mu_k})
\end{equation}

After simplification (\ref{large-delta-1}) becomes
$$2N-4\sum_{m\neq l}\frac{Q^k_l}{Q^k_l-Q^k_m}+\delta_k \bar \nabla (\log \mathfrak{h}_k)(\delta_k Q^k_l)Q^k_l=E. $$
In this section we use $E$ to denote $O(\delta_k^3)+O(\mu_ke^{-\mu_k})$.

According to previous computation (for simplicity we use $m_l$ instead of $m_l^k$ in this section)
\begin{align*}
&\frac{Q^k_l}{Q^k_l-Q^k_j}
=\frac{e^{i\beta_l}(1+m_l)}{e^{i\beta_l}(1-e^{i\beta_{jl}}+m_l-m_je^{i\beta_{jl}})}\\
=&\frac{1}{1-e^{i\beta_{jl}}}-\frac{m_l-m_je^{i\beta_{jl}}}{(1-e^{i\beta_{jl}})^2}+\frac{m_l}{1-e^{i\beta_{jl}}}
+E.
\end{align*}

After simplification we have
$$
\frac{Q^k_l}{Q^k_l-Q^k_j}=\frac{1}{1-e^{i\beta_{jl}}}+\frac{e^{i\beta_{jl}}}{(1-e^{i\beta_{jl}})^2}(m_j-m_l)
+E.
$$

Using (\ref{e-N}) for each $l$, we have

\begin{equation}\label{pi-each-l}
\sum_{j=0,j\neq l}^N\frac{4e^{i\beta_{jl}}}{(1-e^{i\beta_{jl}})^2}(m_j-m_l)
+\delta_k \bar \nabla (\log \mathfrak{h}_k)(\delta_k Q^k_l)e^{i\beta_l}=E
\end{equation}
Using $L_k$ to denote $\nabla \log \mathfrak{h}_k$ and (\ref{trig-i-1}), we write (\ref{pi-each-l}) as
\begin{equation}\label{pi-each-2}
-\sum_{j=0,j\neq l}^Nd_{|j-l|}m_j+Dm_l
+\delta_k \bar L_k(\delta_kQ^k_l)e^{i\beta_l}=E.
\end{equation}
for $l=1,...,N$. Note that in the last term on the left hand side we used $m_l^k=O(\delta_k^2)$, so for this term there is no need to evaluate at $Q_l^k$.

 For $l=0$, using $\beta_0=0$ and $m_0=0$, we have
\begin{equation}\label{pi-each-3}
-\sum_{j\neq 0}^Nd_{j}m_j
+\delta \bar L_k(\delta Q^k_0)=E.
\end{equation}
For the case $N=1$ we have
$$D=1,\quad d_1=1,\quad \beta_{01}=-\pi,\quad \beta_1=\pi,\quad \beta_0=0. $$
Thus (\ref{pi-each-2}) and (\ref{pi-each-3}) are reduced to
\begin{equation}\label{for-N-1}
\left\{\begin{array}{ll}m_1-\delta_k\bar L_k(\delta_kQ^k_1)=E,\\
-m_1+\delta_k\bar L_k(\delta_k Q^k_0)=E.
\end{array}
\right.
\end{equation}
where we have used $m_1^2=O(\delta_k^4)+O(e^{-\mu_k})$.
Adding the two equations in (\ref{for-N-1}), we have
\begin{equation}\label{for-N-1-b}
\delta_k\bar L_k(\delta_k Q^k_0)-\delta_k \bar L_k(\delta_k Q^k_1)=E.
\end{equation}
By $Q_0^k=1+E$ and $Q_1^k=-1+E$ we evaluate $L_k$ as
$$L_k(\delta_kQ^k_0)=\left(\begin{array}{c}
\partial_1\log \mathfrak{h}_k(\delta_k Q^k_0)\\
\partial_2\log \mathfrak{h}_k(\delta_k Q^k_0)
\end{array}
\right)
=\left(\begin{array}{c}
\partial_1 \log \mathfrak{h}_k(0)+\delta_k \partial_{11}\log \mathfrak{h}_k(0)\\
\partial_2 \log \mathfrak{h}_k(0)+\delta_k \partial_{12}\log \mathfrak{h}_k(0)
\end{array}
\right)+O(\delta_k^2),
$$
and
$$L_k(\delta_kQ^k_1)=\left(\begin{array}{c}
\partial_1\log \mathfrak{h}_k(\delta_k Q^k_1)\\
\partial_2\log \mathfrak{h}_k(\delta_k Q^k_1)
\end{array}
\right)
=\left(\begin{array}{c}
\partial_1 \log \mathfrak{h}_k(0)-\delta_k \partial_{11}\log \mathfrak{h}_k(0)\\
\partial_2 \log \mathfrak{h}_k(0)-\delta_k \partial_{12}\log \mathfrak{h}_k(0)
\end{array}
\right)+O(\delta_k^2).
$$
By comparing coefficients we have
\begin{equation}\label{van-n-1}
a_{11}=\delta_k^{-2}E,\quad
a_{12}=\delta_k^{-2}E,
\end{equation}
where $a_{11}=\partial_{11}(\log \mathfrak{h}_k)(0)$, $a_{12}=\partial_{12}(\log \mathfrak{h}_k)(0)$.
Direct computation from
$$\mathfrak{h}_k(x)=h_k(xe^{i\theta_k})$$ gives
\begin{align*}
a_{11}&=\partial_{11}(\log h_k)(0)\cos^2(\theta_k)+2\partial_{12}(\log h_k)(0)\cos\theta_k\sin\theta_k\\
&+\partial_{22}(\log h_k)(0)\sin^2\theta_k
=\partial_{e_ke_k}(\log h_k)(0),
\end{align*}
\begin{align*}
a_{12}&=(\partial_{22}-\partial_{11})(\log h_k)(0)\sin\theta_k\cos\theta_k+\partial_{12}(\log h_k)(0)(\cos^2\theta_k-\sin^2\theta_k)\\
=&\partial_{e_ke^{\perp}_k}(\log h_k)(0).
\end{align*}

\medskip

For the case $N\ge 2$ we take the sum of all equations in (\ref{pi-each-2}), (\ref{pi-each-3}) and obtain
$$\delta_k\sum_{l=0}^N\bar L_k(\delta_k Q_l^k)e^{i\beta_l}=E,$$
This equation is equivalent to
\begin{equation}\label{2nd-va-2}
\sum_{l=0}^N\delta_k(\bar L_k(\delta_k Q^k_l)-\bar L_k(0))e^{i\beta_l}=E,
\end{equation}
where we have used $\sum_{l=0}^Ne^{i\beta_l}=0$.
Using
$$a_{11}=\partial_{11}(\log \mathfrak{h}_k)(0),\quad a_{12}=\partial_{12}(\log \mathfrak{h}_k)(0),\quad a_{22}=\partial_{22}(\log \mathfrak{h}_k)(0),$$
 we see that
\begin{align*}
&\bar L_k(\delta_kQ_l^k)-\bar L_k(0)\\
=&\partial_1 \log \mathfrak{h}_k(\delta_k Q_l^k)-i\partial_2 \log \mathfrak{h}_k(\delta_k Q_l^k)-(\partial_1 \log \mathfrak{h}_k(0)-i\partial_2 \log \mathfrak{h}_k(0))\\
=&(a_{11}\delta_k\cos\beta_l+a_{12}\sin\beta_l)-i(a_{12}\cos\beta_l+a_{22}\sin\beta_l)
\end{align*}
Thus the real part of (\ref{2nd-va-2}) is
$$
\delta_k^2(\sum_{l=0}^Na_{11}\cos^2\beta_l +2a_{12}cos \beta_l \sin \beta_l+a_{22} \sin^2\beta_l^2)=\frac{N+1}2(a_{11}+a_{22})\delta_k^2
$$
and the imaginary part of (\ref{2nd-va-2}) is
$$
\delta_k^2(\sum_{l=0}^N(a_{11}-a_{22})\sin\beta_l\cos\beta_l-a_{12}\cos (2\beta_l))=0,
$$
where we have used
$\sum_{l=0}^Ne^{2i\beta_l}=0. $
Thus (\ref{2nd-va-2}) becomes
$$\sum_{l=0}^N(\bar L_k(\delta_k Q^k_l)-\bar L_k(0))e^{i\beta_l}=\left(\begin{array}{c}
\delta_k^2\frac{N+1}2\Delta (\log \mathfrak{h}_k)(0)\\
0
\end{array}
\right) +E,
$$
and we have
$$\Delta (\log \mathfrak{h}_k)(0)=O(\delta_k^{-2}\mu_k e^{-\mu_k})+O(\delta_k). $$
Since $\Delta (\log h_k)(0)=\Delta (\log \mathfrak{h}_k)(0)$, we obtain the conclusion stated in Theorem \ref{Van-n1}. $\Box$

\section{Appendix A: Simple blowup solutions}

In the section we approximate simple bubbling solutions using global solutions.
Recall that $u_k$ satisfies
$$\Delta u_k+|x|^{2N}h_k(x)e^{u_k}=0,\quad |x|<\tau, $$
$$\max_{\bar B_{\tau}}u_k=\lambda_k\to \infty, $$
$$0 \mbox{ is the only blowup point in } \quad B_{\tau}, $$
and
$$u_k \mbox{ is a constant at }\quad \partial B_{\tau}. $$
Since we talk about simple blowup solutions, we have
\begin{equation}\label{sb-ine}
u_k(x)+2(1+N)\log |x|\le C
\end{equation}
for some $C>0$ independent of $k$.

Let
\begin{equation}\label{a-ep}
\epsilon_k=e^{-\frac{\lambda_k}{2(1+N)}}
\end{equation}
 and
 \begin{equation}\label{a-vk}
v_k(y)=u_k(\epsilon_k y)+2(1+N)\log \epsilon_k,\quad |y|\le 1/\epsilon_k.
\end{equation}
Then clearly $v_k\le 0$ and satisfies
\begin{equation}\label{vk-e}
\Delta v_k+h_k(\epsilon_k y)|y|^{2N}e^{v_k(y)}=0,\quad \mbox{ in }\quad |y|<\tau/\epsilon_k.
\end{equation}

It is easy to see that a subsequence of $v_k$, which is denoted as $v_k$ as well, converges uniformly to $v$ over any fixed compact subset of $\mathbb R^2$.
The limit function $v$, which solves
\begin{equation}\label{v-e}
\Delta v(y)+|y|^{2N}e^{v}=0,\quad \mbox{in }\quad \mathbb R^2, \quad \int_{\mathbb R^2}|y|^{2N}e^v<\infty,
\end{equation}
also satisfies
$$v(y)\le 0 $$
and, by the classification theorem of Prajapat-Tarantello \cite{prajapat}
$$\int_{\mathbb R^2}|y|^{2N}e^{v}=8\pi(1+N). $$
The asymptotic behavior of $v$ is determined by the total integration of $\int_{\mathbb R^2}|y|^{2N}e^{v}$:
$$v(y)=-4(1+N)\log |y|+O(1)\quad \mbox{for}\quad |y|>1. $$

So $v_k\to v$ over any fixed $B_R$ in $\mathbb R^2$. Next we consider the behavior of $v_k$ outside $B_R$.

For $r\in (2\epsilon_k R, \tau/3)$, let
$$\tilde v_k(y)= u_k(ry)+2(1+N)\log r, \quad \frac 12<|y|<2. $$
Then clearly $\tilde v_k$ satisfies
$$\Delta \tilde v_k(y)+|y|^{2N}h_k(ry)e^{\tilde v_k(y)}=0, \quad B_2\setminus B_{1/2}. $$
Let $c_0$ be the bound for $\tilde v_k$: $\tilde v_k\le c_0$ in $B_2\setminus B_{1/2}$ and we set
$g_k=\tilde v_k-c_0-1$ which immediately satisfies $g_k\le -1$ in $B_2\setminus B_{1/2}$. Thus the equation for $g_k$ can be
written as
$$\Delta g_k+\frac{|y|^{2N}h_k(ry)e^{\tilde v_k}}{g_k}g_k(y)=0, \quad \mbox{ in }\quad B_2\setminus B_{1/2}. $$
The coefficient of $g_k$ is clearly bounded. Thus standard Harnack inequality on $\partial B_1$ gives
$$\max_{\partial B_1}(-g_k)\le c_1(c_0) \min_{\partial B_1} (-g_k) $$
where $c_1>1$ only depends on $c_0$.
Going back to $\tilde v_k$ we have
$$\max_{\partial B_1} \tilde v_k\le \frac{1}{c_1}\min_{\partial B_1} \tilde v_k+(c_0+1)(1-\frac{1}{c_1}). $$

For $u_k$ it is
\begin{equation}\label{sphe-uk}
\max_{\partial B_r} u_k\le \frac{1}{c_1}\min_{\partial B_r}u_k-2(1+N)(1-\frac{1}{c_1})\log r+(c_0+1)(1-\frac{1}{c_1}).
\end{equation}

Let $\bar v_k(r)$ be the spherical average of $v_k$ on $\partial B_r$. Then for $r>R$ for some $R$ large,
$$\bar v_k(r)\le (-4(1+N)+\delta_1)\log r +O(1) $$
because
$$\frac{d}{dr}\bar v_k(r)=-\frac{1}{2\pi}\int_{B_r}|y|^{2N}h_k(\epsilon_ky)e^{v_k(y)}dy. $$
From the value of $\bar v_k(R)$ and the estimate above, it is immediate to see that
\begin{equation}\label{v-sphe-a}
\bar v_k(r)\le (-4(N+1)+\epsilon(R))\log r+c, \quad r\ge R.
\end{equation}
For $r\ge R$, the spherical Harnack inequality for $u_k$ gives,
$$\max_{\partial B_r}v_k\le \frac{1}{c_1}\min_{\partial B_r}v_k-2(1+N)(1-\frac{1}{c_1})\log r+c_2 $$
where $c_2=(c_0+1)(1-\frac{1}{c_1})$.
This inequality readily  gives the following estimate of $v_k$:
$$v_k(y)\le (-2(1+N)-\delta)\log |y|,\quad R<|y|<\tau \epsilon_k^{-1} $$
for some $\delta>0$. Then it is easy to use Green's representation formula to obtain
$$v_k(y)\le -4(1+N)\log (1+|y|)+c,\quad |y|\le \tau/\epsilon_k. $$
The classification theorem of Prajapat-Tarantello \cite{prajapat} gives
$$v(y)=\log \frac{\Lambda }{(1+\frac{\Lambda}{8(1+N)^2}|y^{N+1}-\xi|^2)^2}, $$
where parameters $\Lambda>0$ and $\xi\in \mathbb C$.
By the argument in Lin-Wei-Zhang \cite{lwz-jems} there is a perturbation of $\Lambda_k\to \Lambda$ and $\xi_k\to \xi$ such that
\begin{equation}\label{Vk-def}
V_k(y):=\log \frac{\Lambda_k }{(1+\frac{\Lambda_k}{8(1+N)^2}|y^{N+1}-\xi_k|^2)^2},
\end{equation}
satisfies, for carefully chosen $p_1,p_2,p_3$, that
\begin{equation}\label{vk-better}
v_k(p_l)=V_k(p_l),\quad l=1,2,3.
\end{equation}

The idea of the proof in \cite{lwz-jems} for this case is the following: Choose $1<<|p_1|<<|p_2|<<|p_3|$ such that the following matrix invertible:
\begin{equation}\label{ap-ma}
M=\left(\begin{array}{ccc}
\frac{\partial v}{\partial \Lambda}(p_1) & \frac{\partial v}{\partial \Lambda}(p_2) & \frac{\partial v}{\partial \Lambda}(p_3)\\
\frac{\partial v}{\partial \xi_1}(p_1) & \frac{\partial v}{\partial \xi_1}(p_2) & \frac{\partial v}{\partial \xi_1}(p_3)\\
\frac{\partial v}{\partial \xi_2}(p_1) & \frac{\partial v}{\partial \xi_2}(p_2) & \frac{\partial v}{\partial \xi_2}(p_3)
\end{array}
\right)
\end{equation}
where $\xi=\xi_1+i\xi_2$. Thus if a $o(1)$ perturbation is placed on $v$ (to make (\ref{vk-better}) hold), all we need to do is change the parameters $\Lambda$, $\xi$ by a comparable amount. So even though we have a sequence of parameters $\Lambda_k$, $\mu_k$, they are not tending to infinity.

Direct computation shows
\begin{align*}
\frac{\partial v}{\partial \Lambda}=-\frac{1}{\Lambda}+\frac{1}{8(N+1)^2}\frac{1}{1+\frac{\Lambda}{8(N+1)^2}|z^{N+1}-\xi|^2},\\
\frac{\partial v}{\partial \xi}=\frac{\Lambda}{4(1+N)^2}\frac{\bar z^{N+1}-\bar \xi}{1+\frac{\Lambda}{8(1+N)^2}|z^{N+1}-\xi|^2},\nonumber \\
\frac{\partial v}{\partial \bar \xi}=\frac{\Lambda}{4(1+N)^2}\frac{ z^{N+1}-\xi}{1+\frac{\Lambda}{8(1+N)^2}|z^{N+1}-\xi|^2},\nonumber
\end{align*}

Obviously $M$ is invertible if and only if the following matrix is invertible:
$$M_1=\left(\begin{array}{ccc}
\frac{\partial v}{\partial \Lambda}(p_1) & \frac{\partial v}{\partial \Lambda}(p_2) & \frac{\partial v}{\partial \Lambda}(p_3)\\
\frac{\partial v}{\partial \xi}(p_1) & \frac{\partial v}{\partial \xi}(p_2) & \frac{\partial v}{\partial \xi}(p_3)\\
\frac{\partial v}{\partial \bar \xi}(p_1) & \frac{\partial v}{\partial \bar \xi}(p_2) & \frac{\partial v}{\partial \bar \xi}(p_3)
\end{array}
\right)$$
we choose $p_l$ to be
$$p_l=s^{1+\epsilon l}  e^{\sqrt{-1}\theta_l},\quad l=1,2,3 $$
where $s\ge \ge 1\ge \ge \epsilon>0 $ are constants to be determined later. Using crude expansion, we can write $M_1$ as
$$
\left(\begin{array}{ccc}
-\frac{1}{\Lambda}+O(|p_1|^{-2N-2}) & -\frac{1}{\Lambda}+O(|p_2|^{-2N-2}) & -\frac{1}{\Lambda}+O(|p_3|^{-2N-2}) \\
\frac{2\Lambda}{c_0}\frac{e^{-(N+1)i\theta_1}}{|p_1|^{N+1}}+O(\frac{1}{|p_1|^{2N+2}}) & \frac{2\Lambda}{c_0}\frac{e^{-(N+1)i\theta_2}}{|p_2|^{N+1}}+O(\frac{1}{|p_2|^{2N+2}}) &
\frac{2\Lambda}{c_0}\frac{e^{-(N+1)i\theta_3}}{|p_3|^{N+1}}+O(\frac{1}{|p_3|^{2N+2}}) \\
\frac{2\Lambda}{c_0}\frac{e^{(N+1)i\theta_1}}{|p_1|^{N+1}}+O(\frac{1}{|p_1|^{2N+2}}) & \frac{2\Lambda}{c_0}\frac{e^{(N+1)i\theta_2}}{|p_2|^{N+1}}+O(\frac{1}{|p_2|^{2N+2}}) &
\frac{2\Lambda}{c_0}\frac{e^{(N+1)i\theta_3}}{|p_3|^{N+1}}+O(\frac{1}{|p_3|^{2N+2}})
\end{array}
\right)
$$
where for simplicity we use $c_0=8(1+N)^2$.
Multiplying the first row by $-\Lambda$, the second row and the third row by $\frac{c_0}{2\Lambda}|p_3|^{2N+2}$, we change $M_1$ to

$$\left(\begin{array}{ccc}
1+O(s^{-N-1}) & 1+O(s^{-N-1}) & 1+O(s^{-N-1}) \\
s^{2(N+1)\epsilon}e^{-i(N+1)\theta_1}+S_1 & s^{(N+1)\epsilon}e^{-i(N+1)\theta_2}+S_2 &
e^{-i(N+1)\theta_3}+S_3 \\
s^{2(N+1)\epsilon}e^{i(N+1)\theta_1}+S_1 & s^{(N+1)\epsilon}e^{i(N+1)\theta_2}+S_2 &
e^{i(N+1)\theta_3}+S_3
\end{array}
\right)
$$
where
$$S_1=O(s^{-(N+1)(1-\epsilon)})\quad S_2=O(s^{-(N+1)(1+\epsilon)})\quad S_3=O(s^{-(N+1)(1+3\epsilon)}). $$
By choosing $s>>1$ and $0<\epsilon<<1$ it is easy to see that the determinant of the matrix above is not zero if any only if the following matrix is non-zero:
$$M_2:=\left(\begin{array}{ccc}
1 & 1 & 1 \\
s^{2(N+1)\epsilon}e^{-i(N+1)\theta_1} & s^{(N+1)\epsilon}e^{-i(N+1)\theta_2} &
e^{-i(N+1)\theta_3} \\
s^{2(N+1)\epsilon}e^{i(N+1)\theta_1} & s^{(N+1)\epsilon}e^{i(N+1)\theta_2} &
e^{i(N+1)\theta_3}
\end{array}
\right)
$$
The determinant of $M_2$ is
$$det(M_2)=2i \sin((N+1)(\theta_2-\theta_1)) s^{3(N+1)\epsilon}+O(s^{(2N+2)\epsilon}).$$
 Thus by choosing $s$ large, $\epsilon$ small and $\theta_1,\theta_2$ appropriately we can make $M_2$, as well as $M$, invertible.

 Let
$$w_k(y)=v_k(y)-V_k(y),\quad |y|\le \tau/\epsilon_k. $$
The equation for $w_k$ is
\begin{equation}\label{e-wk}
\Delta w_k+|y|^{2N}e^{\xi_k}w_k=-|y|^{2N}(\sum_{t=1}^2\epsilon_k\partial_t h_k(0)y^t+O(\epsilon_k^2|y|^2))e^{v_k}
\end{equation}
for $y\in B(0,\tau\epsilon_k^{-1})$. In addition, we know that $w_k(p_t)=0$ for $t=1,2,3$ and $w_k(y)=O(1)$ and
the oscillation of $w_k$ on $\partial B(0,\tau\epsilon_k^{-1})$ is $O(\epsilon_k^{N+1})$.

Our next step is to improve the estimate of $w_k$.
From the Green's representation formula for $w_k$ we have
\begin{align*}
w_k(y)=\int_{\Omega_k}G(y,\eta)|\eta |^{2N}(e^{\xi_i}w_k(\eta)+\epsilon_k\sum_t\partial_t h_k(0)y^t\\
+O(\epsilon_k^2|\eta|^2)e^{v_k(\eta)})dy
+w_k|_{\partial \Omega_k}+O(\epsilon_k^{N+1}).
\end{align*}
where $\Omega_k=B(0,\tau\epsilon_k^{-1})$. Using crude estimate of $w_k$ we rewrite the above as
\begin{equation}\label{wk-e1}
w_k(y)=\int_{\Omega_k}G(y,\eta)O(\epsilon_k)(1+| \eta |)^{-3-2N}dy+w_k|_{\partial \Omega_k}+O(\epsilon_k^{N+1}).
\end{equation}
Since $w_k(p_1)=0$. Evaluating the above at $p_1$ we have
\begin{equation}\label{wk-e2}
0=\int_{\Omega_k}G(p_1,\eta)O(\epsilon_k)(1+|\eta |)^{-3-2N}d\eta+w_k|_{\partial \Omega_k}+O(\epsilon_k^{N+1}).
\end{equation}
The difference of (\ref{wk-e1}) and (\ref{wk-e2}) gives
\begin{equation}\label{tem-wk}
w_k(y)=\int_{\Omega_k}(G(y,\eta)-G(p_1,\eta))O(\epsilon_k)(1+| \eta |)^{-3-2N}dy
+O(\epsilon_k^{N+1}).
\end{equation}
Then the goal is to prove
\begin{equation}\label{first-o}
w_k(y)=O(\epsilon_k)\log (2+|y|).
\end{equation}

Our argument is by contradiction. Suppose
$$\Lambda_k:=\max_{y\in \bar\Omega_k}\frac{|w_k(y)|}{\epsilon_k\log (2+|y|)}\to \infty$$
and $\Lambda_k$ is attained at $y_k$ (which could appear on $\partial \Omega_k$). We set
$$\hat w_k(y)=\frac{w_k(y)}{\Lambda_k\epsilon_k\log (2+|y_k|)}, $$
which is obviously a solution of
$$\Delta \hat w_k(y)+2|y|^{2N}e^{\xi_k}\hat w_k=\frac{O(1+|y|)^{-4}}{\tilde \Lambda_k \log (2+|y_k|)}$$
and the definition of $\hat w_k$ implies
$$|\hat w_k(y)|\le \frac{\log (2+|y|)}{\log (2+|y_k|)}. $$

Green's representation formula for $\hat w_k(y_k)$ gives
\begin{align*}
&\pm 1=\hat w_k(y_k)=\hat w_k(y)-\hat w_k(p_1)\\
&=\int_{\Omega_k}(G(y_k,\eta)-G(p_1,\eta))(2h_1^k(\epsilon_k\eta)|\eta |^2 e^{\xi_k}\hat w_k(\eta)+\frac{O(1)(1+|\eta|)^{-4}}{\tilde \Lambda_k\log (2+|y_k|)})d\eta\\
&+\frac{O(\epsilon_k)}{\tilde \Lambda_k\log (2+|y_k|)}.
\end{align*}
where we have used $\hat w_k(p_1)=0$ and $\hat w_k=C+O(\epsilon_k^2)$ on $\partial \Omega_k$.

The estimate of the Green's function $G_k$ on $\Omega_k$ is (see \cite{lin-zhang-jfa} for detail)
$$|G(y_k,\eta)-G(p_1,\eta)|\le \left\{\begin{array}{ll}
C(\log |\eta |+\log |y_k|),\quad \eta\in \Sigma_1,\\
C(\log |y|+|\log |y-\eta||,\quad \eta \in \Sigma_2,\\
C|y|/|\eta |\quad \eta\in \Sigma_3.
\end{array}
\right.
$$
where
\begin{align*}
&\Sigma_1=\{\eta\in \Omega_k;\quad |\eta|<|y|/2\quad \}\\
&\Sigma_2=\{\eta\in \Omega_k;\quad |\eta -y|<|y|/2,\quad \}\\
&\Sigma_3=\Omega_k\setminus (\Sigma_1\cup \Sigma_2).
\end{align*}
If $y_k\to y^*$, $\hat w_k$ converges to a solution of
$$\Delta \phi+|y|^{2N}e^{U}\phi=0,\quad \mathbb R^2, $$
with mild growth:
$$|\phi(y)|\le C\log (2+|y|). $$
By the non-degeneracy of the linearized equation,
$$\phi(y)=c_1 \frac{\partial U}{\partial \Lambda}(y)+c_2\frac{\partial U}{\partial \xi_1}(y)+c_3\frac{\partial U}{\partial \xi_2}(y).$$
Using $\phi(p_i)=0$ for $i=1,2,3$, we have, by the invertibility of matrix (\ref{ap-ma}), $c_1=c_2=c_3=0$, thus $\phi\equiv 0$, a contradiction to
$\hat w_k(y_k)=\pm 1$.

If $y_k\to \infty$, the evaluation of $\hat w_k(y_k)=o(1)$ can be obtained by elementary estimates, a contradiction to $\pm 1=\hat w_k(y_k)$. Thus (\ref{first-o}) is established. $\Box$

So the conclusion of this section is
\begin{thm}\label{good-ap}
Let $v_k$, $ V_k$, $\epsilon_k$ be defined in (\ref{a-vk}),(\ref{Vk-def}) and (\ref{a-ep}), respectively, then
$$|v_k(y)-\tilde V_k(y)|\le C\epsilon_k\log (2+|y|),\quad |y|\le \tau \epsilon_k^{-1}. $$
\end{thm}

\section{Appendix B: A frequently used approximation result}

In this appendix we cite Gluck's result about sharp estimates of a sequence of bubbling solutions near a blowup point. The notations in this section are independent.

Let $\Omega$ be a bounded open set in $\mathbb R^2$ that contains an open neighborhood of the origin, let $u_k$ be a sequence of solutions to
$$\Delta u_k+V_k(x)e^{u_k}=0,\quad x\in \Omega\subset\subset \mathbb R^2 $$
with $\int_{\Omega}V_ke^{u_k}\le C$ for some $C>0$ independent of $k$. The assumption on $V_k$ is:
$$\frac 1A\le V_k(x)\le A, \quad \|V_k\|_{C^3(\Omega)}\le A,\quad \mbox{where $A>0$ is a constant}. $$
If $0$ is the only blowup point in $\Omega$ and $u_k$ has bounded oscillation on $\partial \Omega_k$,
the expansion of $u_k$ (satisfying $u_k(0)=\max_{\bar \Omega_k}u_k$ ) around $0$ is
\begin{align}\label{gluck-e}
u_k(x)&=\log \frac{e^{u_k(0)}}{(1+\frac{V_k(0)}8e^{u_k(0)}|x-q^k|^2)^2}+\psi_k(x)\\
&-8\frac{\Delta(\log V_k)(0)}{V_k(0)}\epsilon_k^2(\log (2+\epsilon_k^{-1}|x|))^2+O(\epsilon_k^2\log \frac{1}{\epsilon_k}), \nonumber
\end{align}
where $\psi_k$ is a harmonic function in $\Omega$ that satisfies $\psi_k(0)$ and $u_k-\psi_k=constant$ on $\partial \Omega$, $\epsilon_k=e^{-\frac 12u_k(0)}$,
$q^k$ is the local maximum of $u_k-\psi_k$, which satisfies
$$q^k=2\epsilon_k^2\frac{\nabla V_k(0)}{V_k^2(0)}+O(\epsilon_k^3), $$ and at the origin there is a vanishing estimate:
$$|\nabla (\log V_k+\psi_k)(0)|=O(\epsilon_k^2\log \frac{1}{\epsilon_k}). $$

\end{document}